\documentclass{elsarticle}

\usepackage{lineno,hyperref}
\journal{arXiv.org}

\biboptions{numbers,sort&compress}

\usepackage{amsmath,amssymb}
\usepackage{mathrsfs}
\usepackage{graphicx}
\usepackage{enumerate}
\usepackage{extarrows}
\usepackage{xcolor}
\usepackage{pdflscape}
\usepackage[version=4]{mhchem}
\usepackage{multirow}
\usepackage{enumerate}
\usepackage{subfigure}
\usepackage{tikz}
\usepackage{booktabs}
\usetikzlibrary{calc}

\newcommand{\tikzmark}[1]{\tikz[overlay,remember picture] \node (#1) {};}
\newcommand{\DrawBox}[3][]{%
    \tikz[overlay,remember picture]{
    \draw[black,#1]
      ($(#2)+(-0.5em,2.0ex)$) rectangle
      ($(#3)+(0.75em,-0.75ex)$);}
}

\usepackage[margin=3.0cm]{geometry}
\begin{document}

\begin{frontmatter}

\title{Algorithmic asymptotic analysis: extending the arsenal of cancer immunology modeling}

%\author{Dimitris G. Patsatzis}
%\address{School of Applied Mathematics and Physical Sciences, National Technical University of Athens, 15780 Athens, Greece}

\author[mymainaddress]{Dimitrios G. Patsatzis\corref{mycorrespondingauthor}}
\address[mymainaddress]{School of Chemical Engineering, National Technical University of Athens, 15772, Athens, Greece}
\cortext[mycorrespondingauthor]{Corresponding author}
\ead{dpatsatzis@mail.ntua.gr}

\begin{abstract}
The recent advances in cancer immunotherapy boosted the development of tumor-immune system models aiming to provide mechanistic understanding and indicate more efficient treatment regimes.~However, the complexity of such models, their multi-scale dynamics and their overparameterized character renders them inaccessible for wide utilization.~In this work, the dynamics of a fundamental model formulating the interactions of tumor cells with natural killer cells, CD8$^+$ T cells and circulating lymphocytes is examined.~It is first shown that the long-term evolution of the system towards high-tumor or tumor-free equilibria is determined by the dynamics of an initial \emph{explosive stage} of tumor progression.~Focusing on this stage, the algorithmic Computational Singular Perturbation methodology is employed to identify the underlying mechanisms confining the system's evolution towards the equilibrium and the governing slow dynamics along them.~It is shown that these insights are preserved along different tumor-immune system and patient-dependent realizations.~Utilizing the obtained mechanistic understanding, a novel reduced model is constructed in an algorithmic fashion, which accurately predicts the dynamics of the system during the explosive stage and includes half of the parameters of the detailed model.~This present analysis demonstrates the potential of algorithmic asymptotic analysis to simplify the complex, overeparameterized and multi-scale nature of cancer immunology models and to indicate the interactions and cell types to target for more effective treatment development.
\end{abstract}


\end{frontmatter}

%%%%%%%%%%%%%%%%%%%%%%%%%%%%%%%%%%%%%%%%%%%%%%%%%%%%%%%%%%%%%%%%%%%%%%%%%%%%%%%%%%%%%%%%%%%%%%%%%%%%%%%%%%%%%%%%%%%

%\linenumbers
%%%%%%%%%%%%%%%%%%%%%%%%%%%%%%%%%%%%%%%%%%%%%%%%%%%%%%%%%%%%%%%%%%%%%%%%%%%%%%%%
\section{Introduction}

One of the leading causes of death worldwide is cancer, counting millions of new diagnosed cases every year, with the incidence and mortality rates rapidly increasing \citep{bray2018global}.~Cancer treatment mainly involves surgery, chemotherapy and radiation, that are frequently combined with complementary therapies, among which special consideration is given to immunotherapy.~During the last decade, the revolution of cancer immunotherapy indicated new strategies for enabling specific anti-tumour responses, indicatively among them the immune checkpoint inhibitors \citep{pardoll2012blockade}, dendritic cell vaccines \citep{hammerstrom2011cancer} and adoptive T cell transfer \citep{baruch2017adoptive}.~Demonstrating increased efficacy and great promise, various cancer immunotherapy treatment regimes have recently gained approval by the FDA \cite{zheng18,hargadon18}, accompanied by thousands of clinical trials \cite{tang18}.~However, due to the lack of mechanistic understanding of the complex interactions between the immune system and the tumor, critical challenges emerge relating to the low treatment response rates, the accurate prediction of treatment efficacy and the possible adverse effects after treatment \citep{lipson2015antagonists,szeto19,peskov19}.

With the prospect to gain mechanistic understanding, various systems-level modeling approaches have been utilized to explore the interactions of the immune system with cancerous tumor \cite{werner14}, including both deterministic and stochastic ordinary differential equations (ODEs) \cite{de2006mixed,osojnik20,kuznetsov94,den2016re,kirschner98,li16,wilson12,hu2018dynamics,rihan14,xu2013stochastic}, partial differential equations (PDEs) \cite{matzavinos2004mathematical,alTameemi12,lai2017combination}, agent-based models \cite{norton19,owen2011mathematical} and data driven statistical modeling \cite{charoentong17}.~The most widely utilized approach, also employed in this work, is the deterministic ODE modeling of spatially averaged population dynamics models, since it (i) can formulate the numerous biological processes taking place among individual cells, (ii) avoids the complexity introduced by the multiple spatial scales and (iii) is particularly insightful for anti-cancer drug development and treatment optimization \cite{singh17}.~Detailed reviews on the dynamics of the tumor-immune system using ODE modeling is provided in Refs.~\cite{mahlbacher19,wilkie2013review,eftimie2011interactions}. 

Despite the overwhelming interest in tumor-immune system modeling, especially during the last 5 years, many obstacles still remain, owing their nature to traditional drawbacks of mathematical biology modeling pipeline.~The first one relates to the balance between oversimplification and overparameterization that should be preserved so that the model can effectively capture the biological reality of tumor dynamics, but also include a few number of parameters in order to be accessible for wide utilization \citep{peskov19}.~As a result, the frequently utilized models include only a limited number of immune system cell types, nonetheless with many parameters \cite{wilson12,li16,kirschner98,mahlbacher19}.%~Following the model development, the next step is to determine the conditions under which specific tumor phenomena arise, such as tumor remission, escape from immune surveillance, relapse after treatment, etc.~Such a determination is frequently accomplished by bifurcation analysis techniques \citep{de2006mixed,osojnik20,kirschner98,kuznetsov94,den2016re}, which however requires intuition to select the proper parameter; a task that becomes significantly cumbersome for detailed models with many parameters.~
~Following the model development and assuming that the model successfully captures biological reality, its predictive power is usually assessed via sensitivity analysis \cite{li16,hu2018dynamics,wilson12,rihan14}, which despite being systematic, it is computationally expensive and does not distinguish between the fast and slow evolution of the system.~This is of major importance, since the latter approach cannot provide any insight into the underlying mechanisms that control tumor dynamics; e.g., the interactions and the cell types that affect tumor size, the time period over which these interactions occur, etc.~Such insights are particularly significant when accounting for anti-cancer therapies, since they can indicate which interactions and cell types to target for more effective treatment regimes.~Finally and most importantly, irrespectively of their size, their detail on capturing tumor specific behaviors and their predictive ability, the models formulating tumor-immune system interactions share a common ground: their inherent multi-scale nature, which originates from the wide range of timescales over which the immune system and tumor interact.~This fast/slow separation of timescales has been exploited by traditional paper-pencil asymptotic analysis techniques, such as the \emph{Quasi-Steady State Approximation} (QSSA) and \emph{matched asymptotics} \citep{osojnik20,kuznetsov94,li16,louzoun2014mathematical,rihan14}.~However, such techniques are hindered by the complexity and the size of the mathematical model under consideration, so that their utilization is limited to small models.

With the aim to provide a systematic and algorithmic framework to address these aforementioned challenges in tumor-immune system dynamics, here the concepts of the \emph{Geometric Singular Perturbation Theory} (GSPT) \citep{fenichel1979geometric,kaper1999systems} are adopted.%, which can successfully systematize various asymptotic analysis techniques, among which are the QSSA and matched asymptotics techniques.
~According to GSPT, multi-scale systems are confined to evolve along low-dimensional \emph{Slow Invariant Manifolds} (SIMs) that emerge in the phase space as a result of the action of the fast timescales.~The evolution of the system on the SIM is governed by a reduced system that is characterized by the slow timescales \citep{verhulst2006nonlinear,kuehn2015multiple}.~In order to identify the components of the model contributing to the formation of the SIM and to the slow system that governs the flow on it, the \emph{Computational Singular Perturbation} (CSP) methodology \citep{lam1994,lam1989} and its tools are employed.

CSP has been applied to provide mechanistic understanding in a wide variety of multi-scale models, mainly in chemical kinetics of reactive flows \citep{manias2016mechanism,valorani2006automatic}, but also in systems biology, pharmacokinetics \citep{patsatzis2019new,patsatzis2016asymptotic} and brain metabolism \citep{patsatzis2019computational}.~To this day, CSP has never been employed for population dynamics models formulating cell-cell interactions, so that the first and major objective of this work is to extend the application of CSP to cancer immunology, thus introducing the algorithmic tools of GSPT in the field.~For this reason, the model in Ref.~\citep{de2006mixed} is adopted, which is considered one of the most influential and fundamental models on the field \citep{de2009mathematical,bellomo2008foundations,alfonso2017biology} and has been extended for the study of immunotherapy treatment in pancreatic cancer patients \citep{louzoun2014mathematical,li2016mathematical}.~Calibrated against human data, the model in Ref.~\citep{de2006mixed} was proposed to explain the dynamics of tumor cells when interacting with natural killer cells, CD8$^+$ T cells and circulating lymphocytes, accounting for both chemotherapy and immunotherapy treatments.~Considering many interactions between the tumor and the immune system (20 parameters without treatment), this model is particularly suitable to demonstrate the ability of CSP to provide mechanistic insight, thus simplifying the inherent complexity of the model.~The second objective of this work is to identify the mechanisms controlling tumor progression and determine the decisive factors leading to tumor persistence/remission.~On top of these identifications, the capability of constructing a reduced model that accurately reproduces the dynamics of the detailed model, using a fewer number of parameters is examined.~It is crucial to note that these objectives are pursued in a purely algorithmic manner, so that this CSP-based framework can be extended to more complex models, without being hindered by the size and the complexity of the model, or researcher's intuition.

The article is organized as follows.~In Section~\ref{sec:Mod}, the mathematical model is presented and the tumor-free and high-tumor equilibia of the system and their stability are briefly discussed.~Next, in Section~\ref{sec:exp}, the period determining the long-term evolution of the system towards an equilibrium is identified and the multi-scale nature of the model during this period is established.~Following this outcome, in Section~\ref{sec:CSPdyn}, CSP is employed to reveal the underlying mechanisms that drive tumor progression and determine the long-term governing dynamics of the system.~Cross-validation of these findings is performed in a dual manner: (i) by utilizing the mechanistic understanding obtained in order to accurately predict the response of the system in parameter perturbations and (ii) by demonstrating that this dynamical behavior is preserved along different initial conditions and parameter sets, indicating the robustness of the results in different tumor-immune system and patient-dependent realizations.~Finally, a reduced model is constructed and its validity and accessibility is briefly discussed in Section~\ref{sec:redmod}, before concluding with the at-home messages this work and its future perspectives.

%%%%%%%%%%%%%%%%%%%%%%%%%%%%%%%%%%%%%%%%%%%%%%%%%%%%%%%%%%%%%%%%%%%%%%%%%%%%%%%%
\section{The mathematical model}
\label{sec:Mod}

\begin{table}[!b]
\centering
\small
\begin{tabular}{c c l}
\multicolumn{3}{l}{\bf Growth rates} \\		
\hline
R$^1$	&	$a T (1 - b T)$					& tumor cells \\
R$^2$	& 	$e C$						& NK  cells \\
R$^3$	& 	$\alpha$						& circulating lymphocytes \\[6pt]
\multicolumn{3}{l}{\bf Death rates} \\
\hline		
R$^4$	& 	$f N$						& NK  cells \\
R$^5$	& 	$m L$						& CD8$^+$ T cells \\
R$^6$	& 	$\beta C$						& circulating lymphocytes \\[6pt]
\multicolumn{3}{l}{\bf Fractional cell kill rates} \\
\hline
R$^7$	& 	$c N T$						& tumor cells by NK cells \\
R$^8$	& 	$D T$						& tumor cells by CD8$^+$ T cells \\[6pt]
\end{tabular}
\begin{tabular}{c c l}
\multicolumn{3}{l}{\bf Recruitment rates} \\
\hline
R$^9$	& 	$g \dfrac{T^2}{h+T^2} N$			& NK cells \\[6pt]
R$^{10}$	& 	$j \dfrac{D^2 T^2}{k+D^2 T^2} L$	& CD8$^+$ T cells  \\[6pt]
R$^{11}$ 	& 	$r_1 N T$						& CD8$^+$ T cells  \\
R$^{12}$ 	& 	$r_2 C T$						& CD8$^+$ T cells \\[6pt]
\multicolumn{3}{l}{\bf Inactivation rates} \\
\hline
R$^{13}$ 	& 	$p N T$						& NK cells \\
R$^{14}$ 	& 	$q L T$						& CD8$^+$ T cells \\
R$^{15}$ 	& 	$u N L^2$						& CD8$^+$ T cells \\
& & \\
\end{tabular}
\caption{The growth, death, fractional cell kill, recruitment and inactivation processes accounted for in the model and the expressions of their reaction rates.}
\label{tb:RR}
\end{table}

The dynamics of the tumor and immune system interactions is examined here in the context of spatially averaged cell populations.~One of the most fundamental mathematical models in cancer-immunology \cite{de2006mixed} is adopted, which formulates the interactions of tumor cells with three immune system cell types: the natural killer (NK) cells, CD8$^+$ T cells and circulating lymphocytes.~Growth, death, fractional tumor cell kill, recruitment and inactivation processes are accounted for in the model, which formulates the evolution of the four cell populations according to the following system of ODEs:
\begin{align}
\dfrac{dT}{dt} 	& = R^1 - R^7 - R^8 &	\dfrac{dN}{dt} 	& = R^2 - R^4  + R^9 - R^{13}  \nonumber \\
\dfrac{dL}{dt} 	& = - R^5 + R^{10} + R^{11} + R^{12} - R^{14} - R^{15}	&  	\dfrac{dC}{dt} 	& = R^3 - R^6  \label{eq:VecF} 
\end{align}
where $T$, $N$, $L$ and $C$ denote the cell population of tumor cells, NK cells, CD8$^+$ T cells and circulating lymphocytes, respectively and $R^k$ denotes the reaction rate of the $k$-th process.~The processes and the expressions of their related reaction rates as enlisted in Table~\ref{tb:RR}, where $D=d \left( (L/T)^l \right) / \left( s+(L/T)^l \right)$. The model in Eq.~\eqref{eq:VecF} is calibrated against human data and can accommodate chemotherapy and immunotherapy drug intervention terms, which are not considered for the present analysis.~Details on the model development and calibration can be found in Refs.~\citep{de2006mixed,de2009mathematical}.

Considering only feasible (non-negative real cell populations) solutions, the system in Eq.~\eqref{eq:VecF} exhibits one tumor-free equilibrium (TFE) and possibly multiple high-tumor equilbria (HTE).~As shown in \ref{app:Eq+St}, the TFE can be analytically calculated as $E_{0}=(T_0,N_0,L_0,C_0)=(0,\alpha e/\beta f,0,\alpha/\beta)$ and it is stable if and only if $(a-d) \beta f< \alpha c e$.~The number of HTE and their stability is numerically determined, depending on the values of the adopted parameter set.  

\begin{figure}[!h]
\centering 
\subfigure[Tumor cells evolution in time]{\includegraphics[scale=0.23]{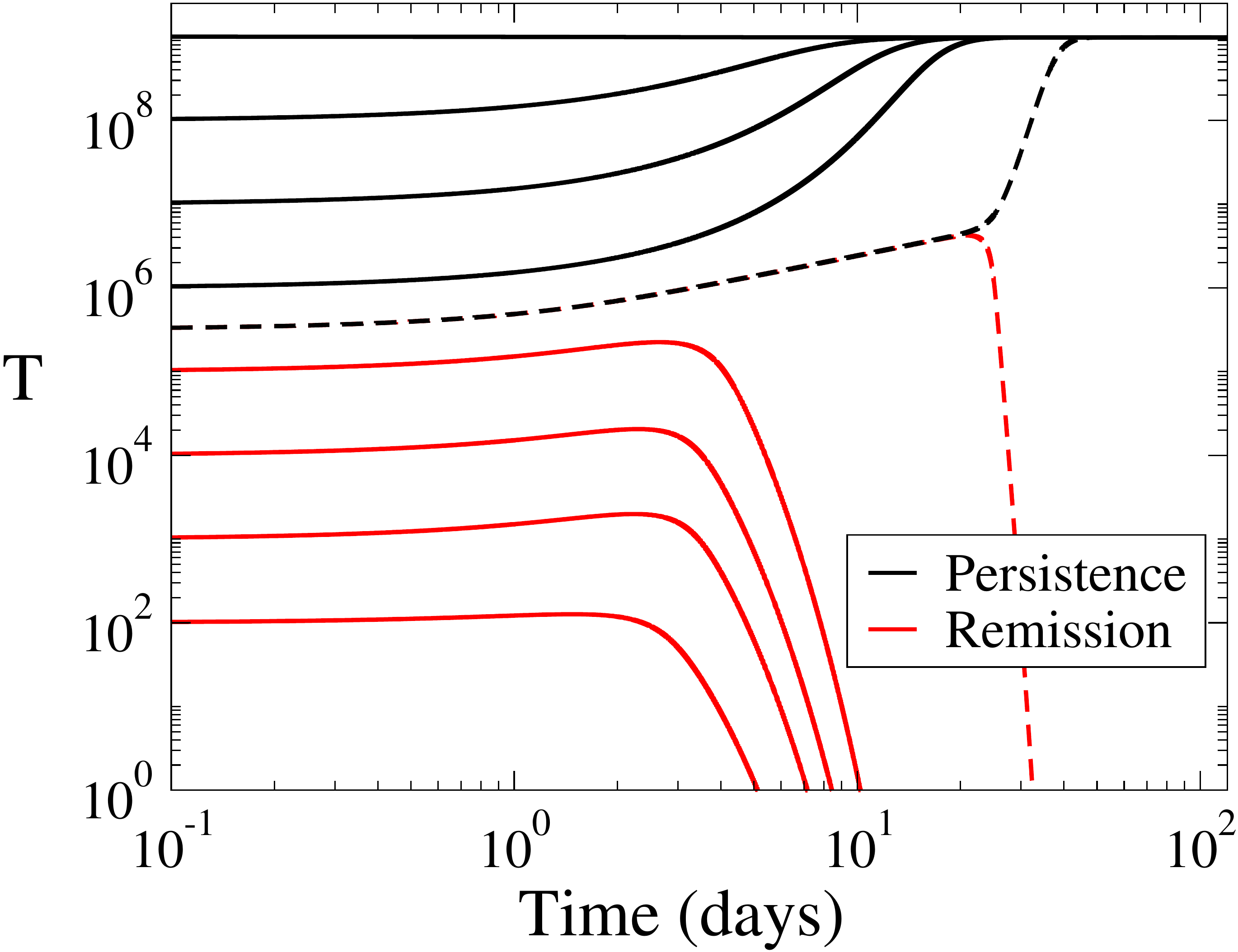}}  \hspace{1.5cm} \subfigure[Bifurcation diagram]{\includegraphics[scale=0.23]{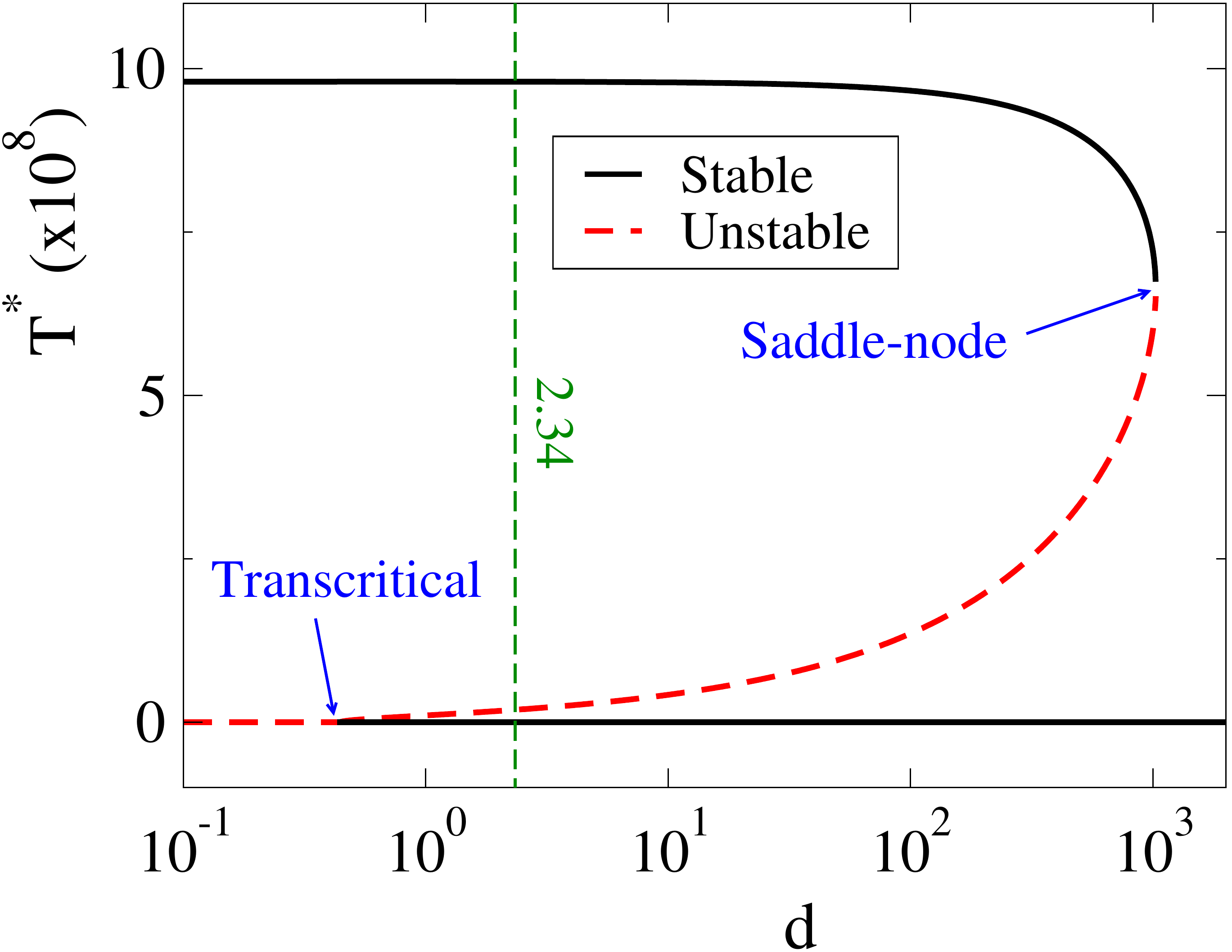}} 
\caption{(a) Response to multiple initial conditions: varying $T(0)$ per order of magnitude (solid lines) and constant $N(0)=10^{+3}$, $L(0)=10^{+1}$, $C(0)=6 \times 10^{+8}$; 1 cell difference (black dashed line $T(0)=319393$ and red dashed line $T(0)=319392$).~(b) Equilibria and their stability with respect to variations of $d$; vertical dashed green line corresponding to the parameter set of patient 9.}
\label{fig:EqST}
\end{figure}

In order to provide an indicative behavior of the system, the parameter set of patient 9 in Ref.~\cite{de2006mixed} was adopted; the parameter values are depicted in \ref{app:ParSet}.~The existence of a stable TFE and a stable HTE, attracting multiple neighboring solutions, is highlighted in Fig.~\ref{fig:EqST}a.~In particular, it is shown that when the initial tumor size is small enough, the system is driven towards the TFE, which is stable due to the validity of the stability condition $(a-d) \beta f< \alpha c e$.~This behavior indicates the ability of the immune system to suppress the tumor.~However, when the initial tumor size is higher than a specific threshold value, the system is driven towards the HTE, which is stable as well, indicating the ability of the tumor to escape immune surveillance.~This dynamical behavior is validated by the bifurcation diagram in Fig.~\ref{fig:EqST}b, in which the existence of the stable TFE and HTE is indicated for the corresponding to patient 9 parameter value of $d=2.34$.~In addition, it is shown that another HTE arises, which however is unstable and thus, does not attract the solutions of the system.

The bifurcation diagram in Fig.~\ref{fig:EqST}b additionally denotes the existence of two bifurcations.~The first is a transcritical bifurcation, located at the point where $(a-d) \beta f = \alpha c e$.~Before this point, the TFE is unstable, so that any initial tumor size is attracted to the HTE, say $E_1$, which is stable.~However, after this point the TFE becomes stable and a new unstable HTE, say $E_2$, emerges.~This new branch sets a threshold, so that (i) any smaller than $E_2$ tumor size is eventually driven towards the TFE, expressing essentially the ability of the immune system to suppress the tumor, while (ii) any larger than $E_2$ tumor size is driven towards the stable HTE $E_1$, expressing essentially the incapability of the immune system to suppress the tumor, which escapes immune surveillance.~The second bifurcation is a saddle-node, after which both HTE $E_1$ and $E_2$ disappear, so that any initial tumor is attracted to the TFE, irrespectively of its size.~Such a behavior is not biologically consistent, a feature that is effectively captured by the extremely high parameter values of $d$.

This insight is very beneficial to understand tumor behavior and its interaction with the immune system.~However, it is informative only for the long-term evolution of tumor behavior, without revealing the underlying mechanisms that drive the system towards it.~In the following Sections, the origin of these mechanisms is examined.

%%%%%%%%%%%%%%%%%%%%%%%%%%%%%%%%%%%%%%%%%%%%%%%%%%%%%%%%%%%%%%%%%%%%%%%%%%%%%%%%
\section{Tumor remission/persistence driving force: the explosive stage}
\label{sec:exp}

In order to examine the decisive factors driving the system towards the HTE or the TFE, two indicative scenarios were considered.~The first scenario represents a ``weakened" immune system, in which the initial tumor persists and escapes immune surveillance and the system is led towards its stable HTE.~The initial cell populations of the tumor persistence case are $N(0)=10^3$, $L(0)=10$ and $C(0)=6\times10^8$.~The second scenario represents a ``healthy" immune system that is able to suppress the initial tumor and the system is led to its TFE.~The initial cell populations of the tumor remission case are $N(0)=10^5$, $L(0)=10^2$ and $C(0)=6\times10^{10}$.~In both the tumor persistence and remission cases, the initial tumor population is $T(0)=10^6$ adopting the initial conditions in \cite{de2006mixed}.

\begin{figure}[!h]
\centering
\includegraphics[scale=0.23]{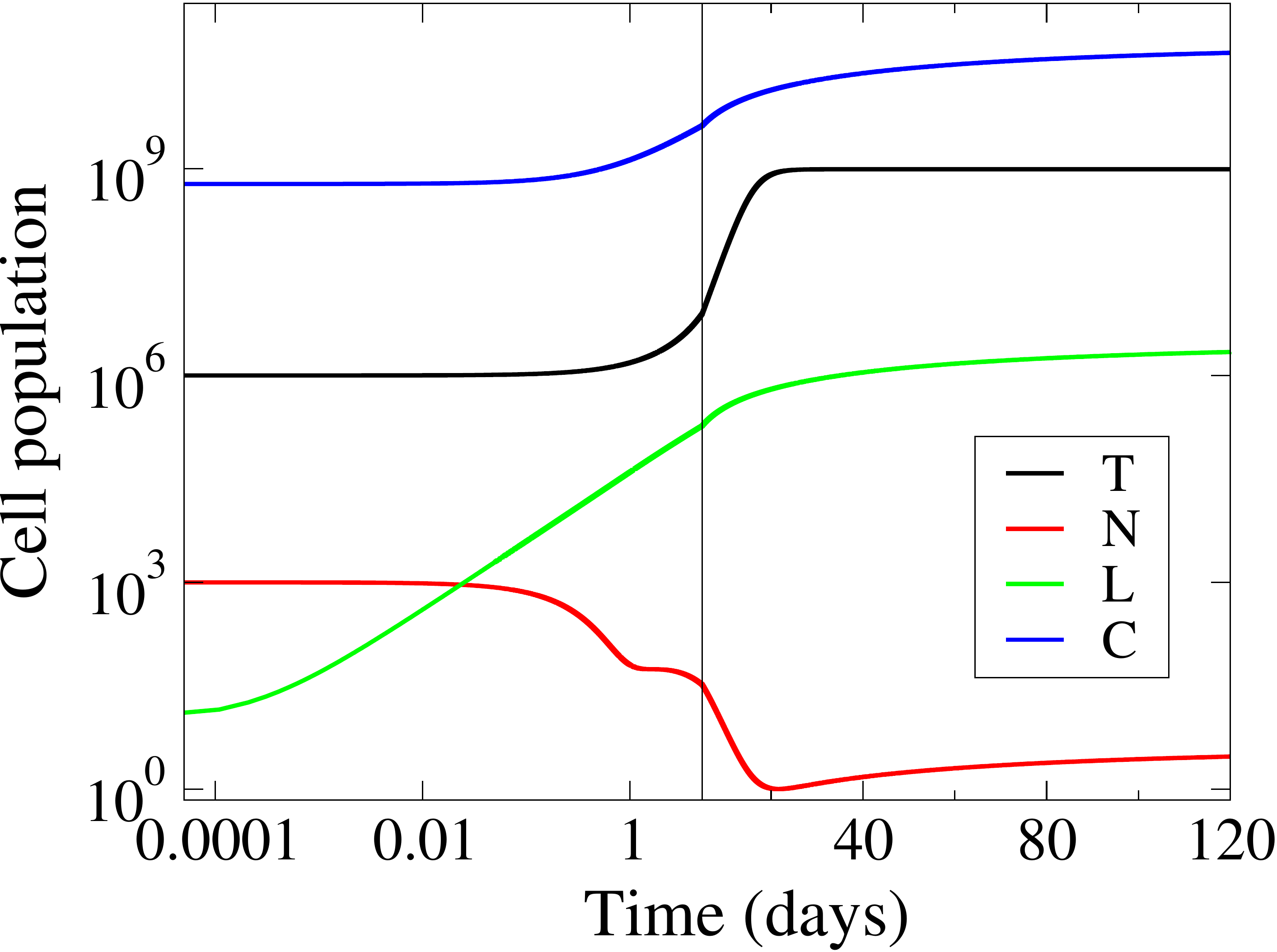} \hspace{1.5cm} \includegraphics[scale=0.23]{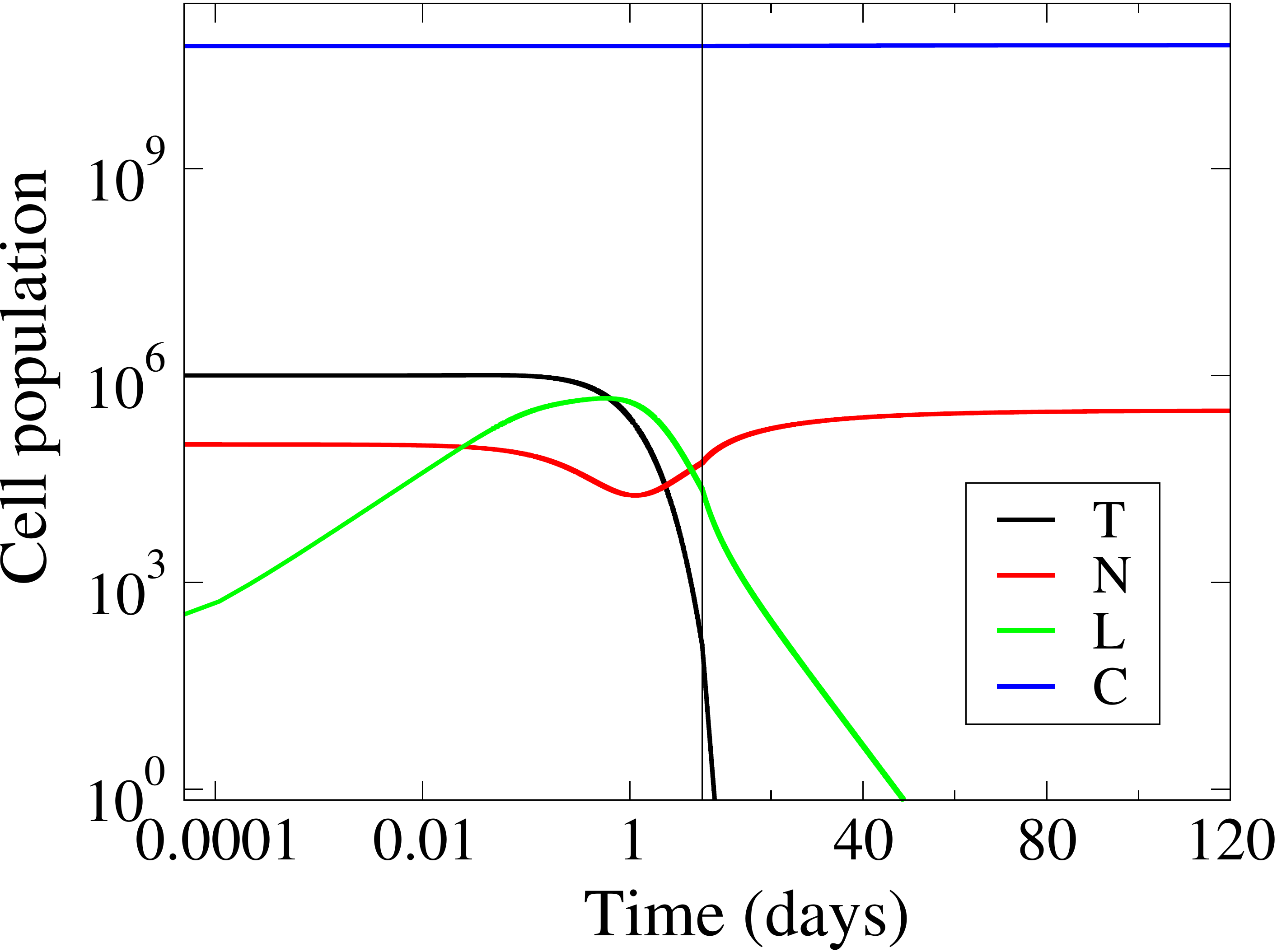} \\
\subfigure[Tumor persistence case]{\includegraphics[scale=0.23]{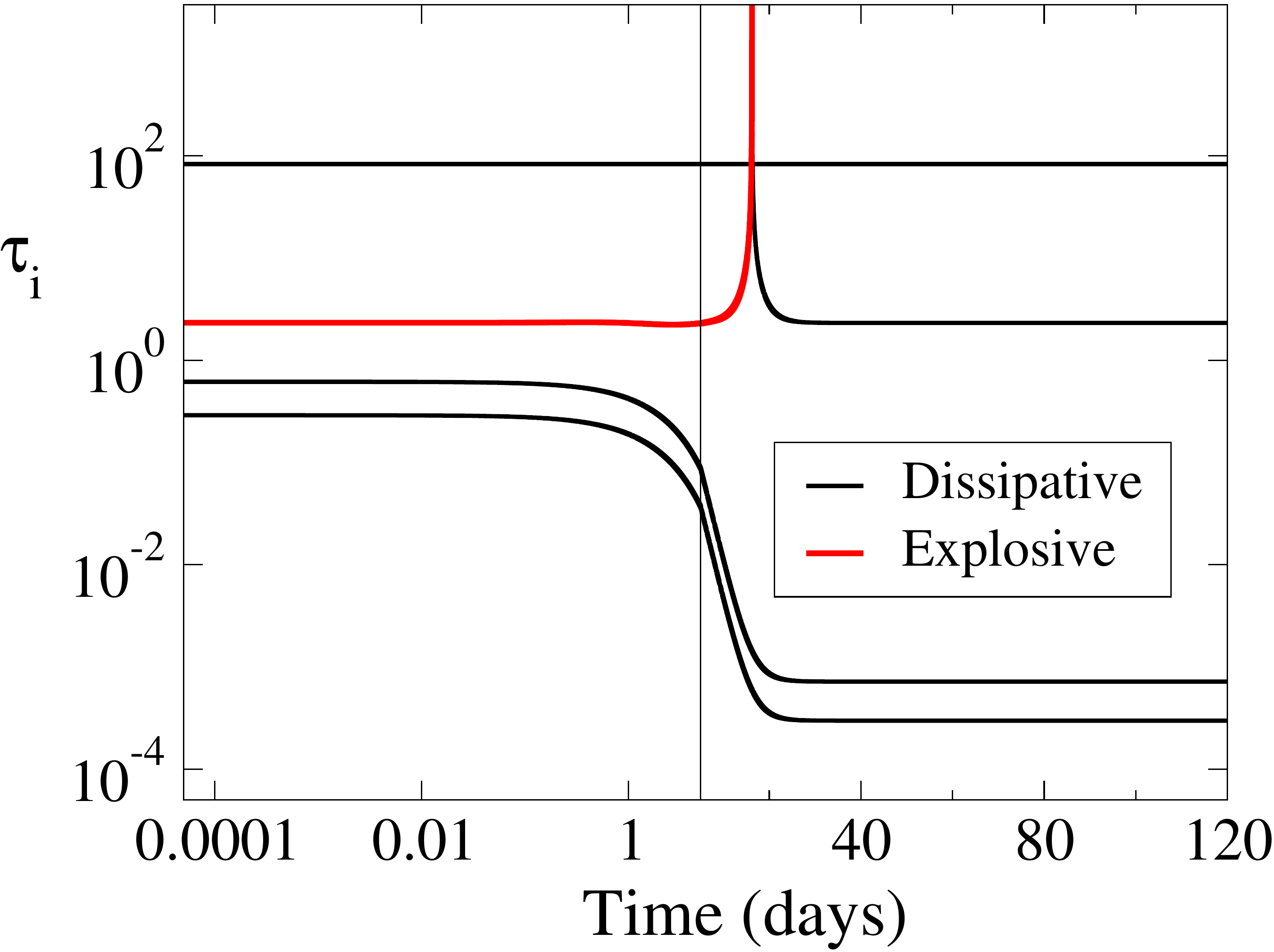}} \hspace{1.5cm} \subfigure[Tumor remission case]{\includegraphics[scale=0.23]{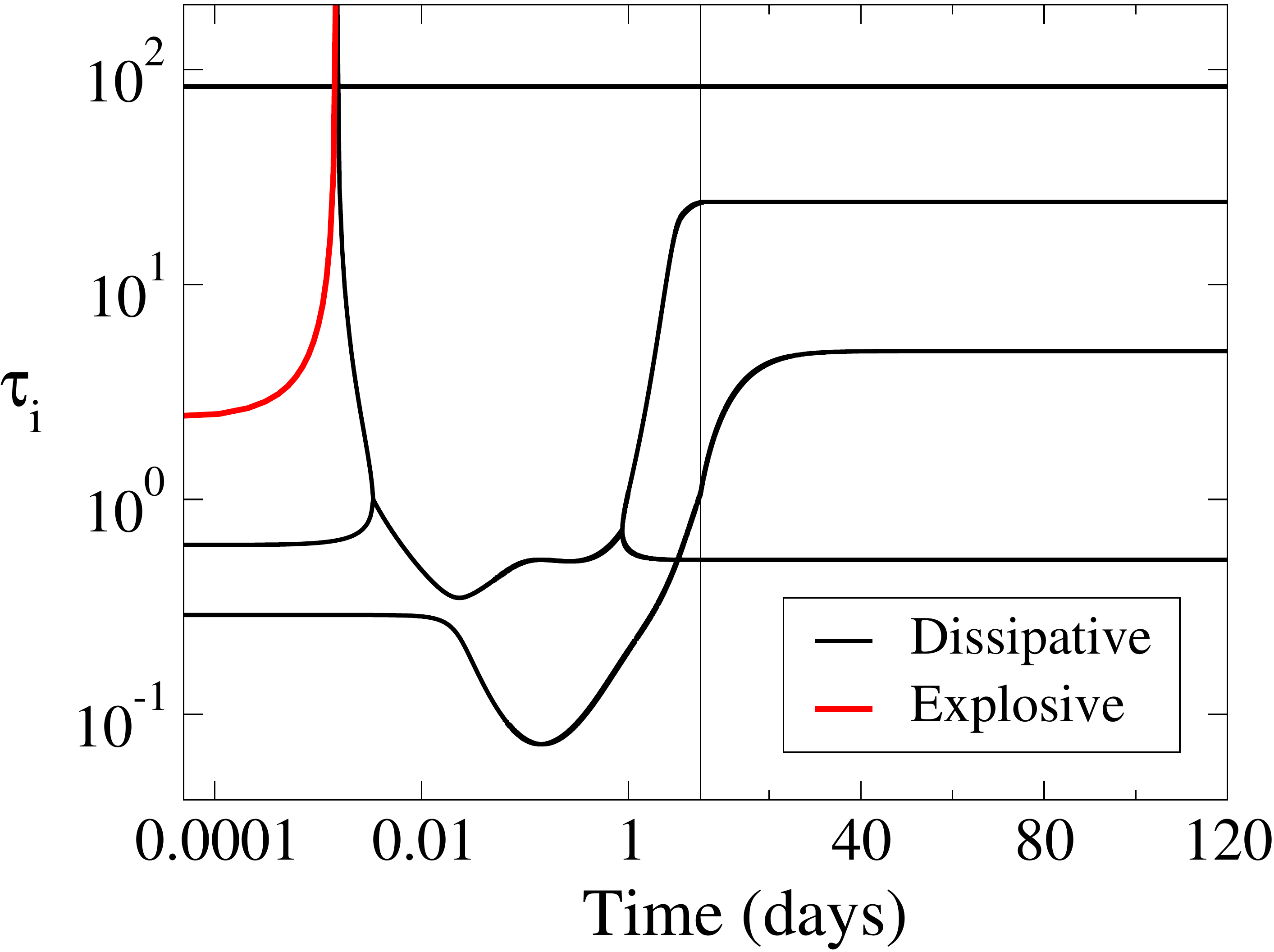}}
\caption{Cell populations (top) and timescales (bottom) of tumor remission and tumor persistence cases.~Logarithmic scale up to the fist 5 days and in linear scale thereafter.}
\label{fig:RPSol+Tmscls}
\end{figure}

The profiles of the cell populations are displayed in the top panels of Fig.~\ref{fig:RPSol+Tmscls}, where it is shown that the initial period up to the first day is qualitatively similar in both cases.~Indeed, despite the differences in the initial conditions, the population of (i) tumor cells, $T$, and circulating lymphocytes, $C$, remain almost constant, (ii) NK cells, $N$, initially remains constant and then decreases, and (iii) CD8$^+$ T cells, $L$, increases quickly.~However, after this initial transient period, the cell population profiles vary significantly.~In particular, in the tumor persistence case, after the first day $L$ continues to increase and $N$ decreases significantly, after attaining a plateau.~The decrease of $N$ results to an increase of $T$ during first 20 days, which leads the system towards its HTE $E_{1}=(9.8 \times 10^{8}, 3.87,2.86 \times 10^6,6.25 \times 10^{10})$, as shown in Fig.~\ref{fig:RPSol+Tmscls}a.~In contrast, in the tumor remission case, $L$ reaches a maximum and then decreases and $N$ reaches at the same time to a minimum value and then increases.~As a result, $T$ faces a rapid decrease up to the first 5 days, which leads the system gradually towards its TFE $E_{0}=(0,3.15 \times 10^{5},0,6.25 \times 10^{10})$, as shown in Fig.~\ref{fig:RPSol+Tmscls}b.

Since in both tumor persistence and remission cases the system in Eq.~\eqref{eq:VecF} is led to its stable HTE and TFE respectively, the behavior of the timescales of the system was examined, as shown in the bottom panels of Fig.~\ref{fig:RPSol+Tmscls}.~In both cases, the system exhibits 4 timescales, which vary from $O(10^{-1})$ to $O(10^2)$ from the beginning of the process.~In the tumor remission case, this range is preserved throughout the process, while in the tumor persistence case, after the first 5 days the timescales extend in a much wider range from $O(10^{-3})$ to $O(10^2)$.~Thus, the multi-scale character of the model is more pronounced in the tumor persistence case.~The most significant difference though, is the nature of the timescales; a dissipative/explosive timescale tends to drive the system towards to/away from equilibrium (negative/positive real part of the corresponding eigenvalue).~In both tumor persistence and remission cases, the system exhibits 3 dissipative and 1 explosive timescales, the later degenerating to a dissipative timescale later on, as well.~However, the period in which the explosive timescale exists is significantly longer in the tumor persistence case ($16.2~d$), in comparison to the tumor remission case ($0.002~d$).

\begin{figure}[!h]
\centering
\subfigure[Solution from Fig.~\ref{fig:RPSol+Tmscls}a,b]{\includegraphics[scale=0.23]{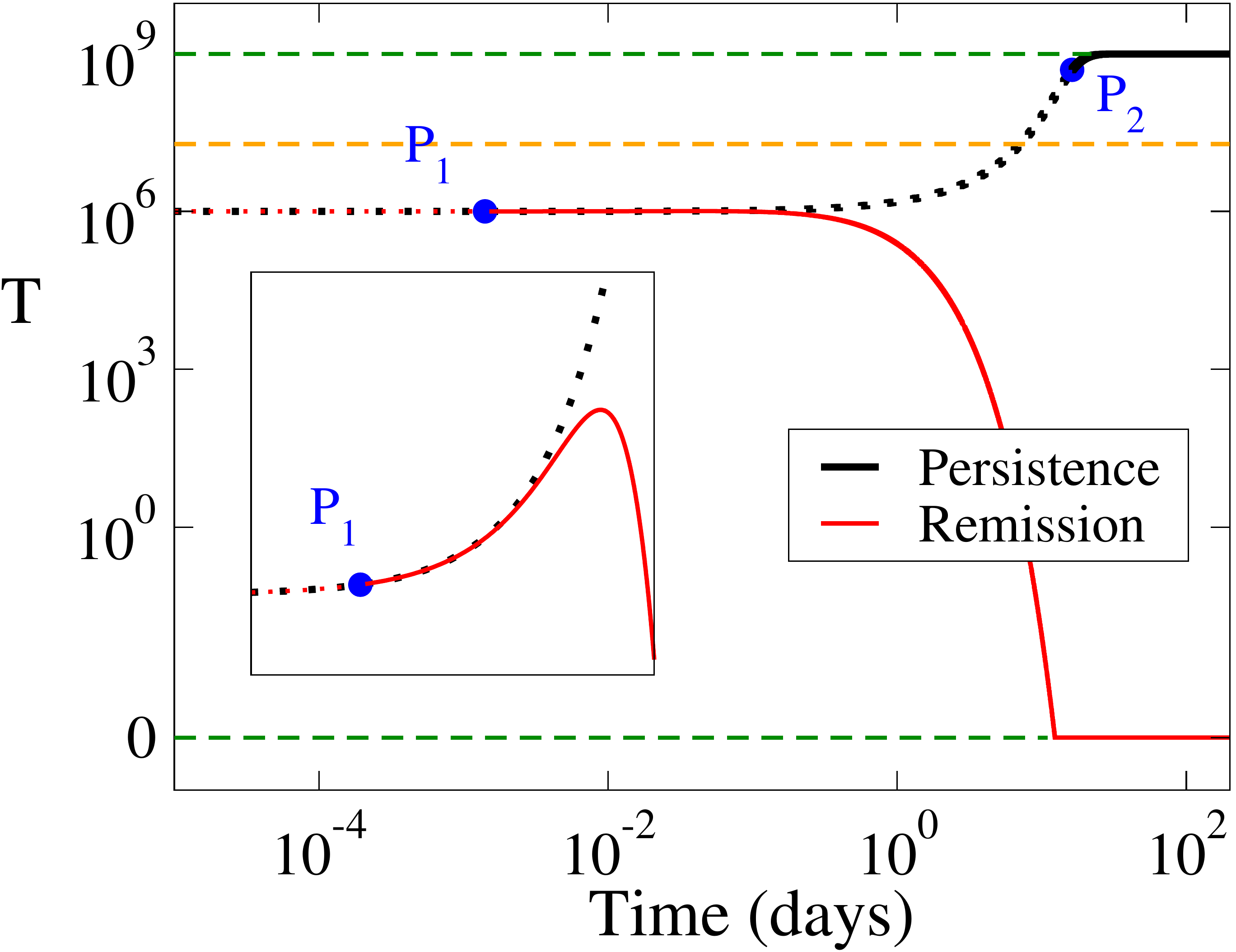}} \hspace{1.5cm} \subfigure[1 initial tumor cell difference]{\includegraphics[scale=0.23]{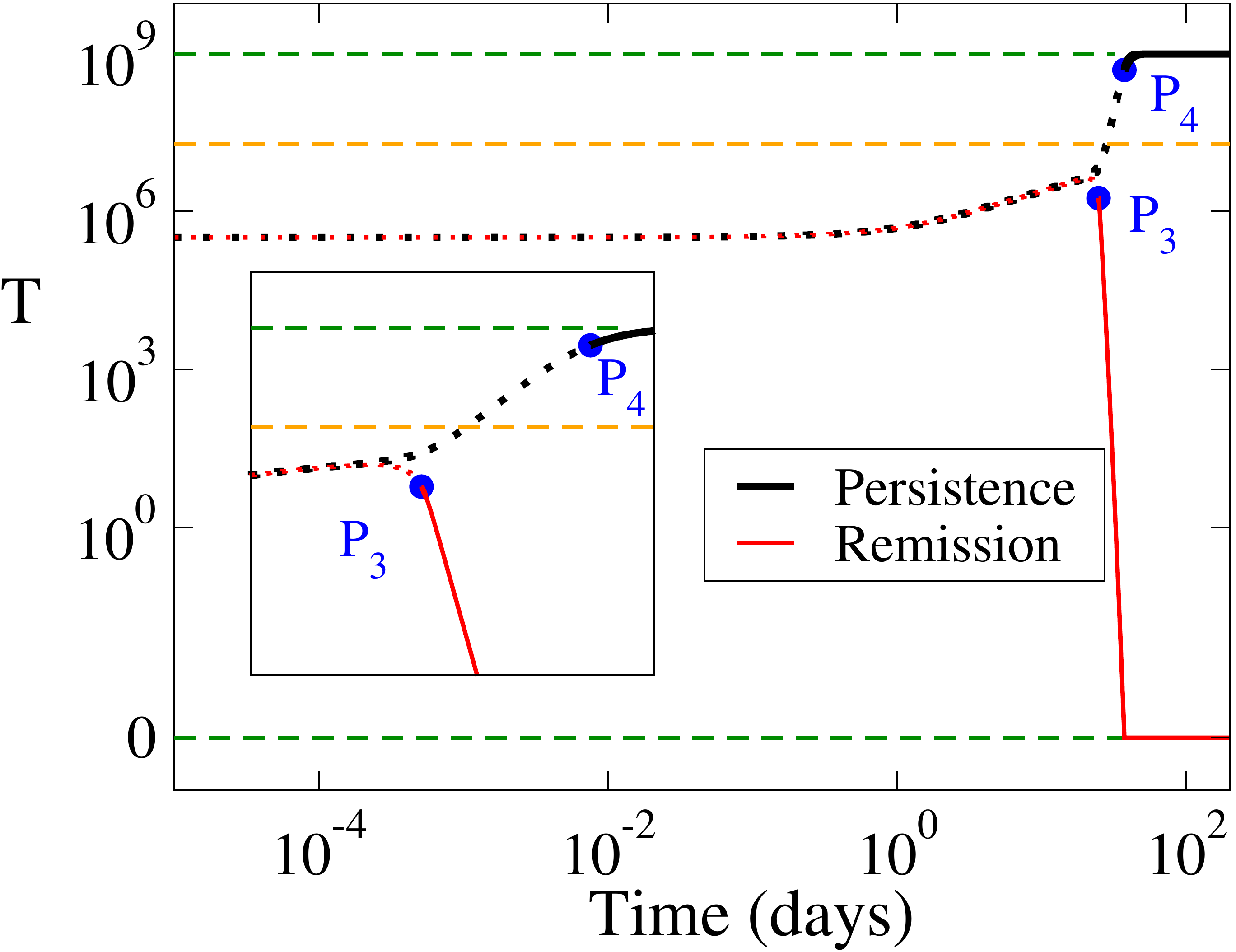}}
\caption{Tumor cell population profiles in persistence (black) and remission (red) cases; dotted during the explosive stage and solid after it, separated by the points P$_1$-P$_4$.~The horizontal dashed lines indicate the TFE and HTE of the system; green for the stable ones and orange for the unstable one.}
\label{fig:RPSolComp}
\end{figure}

In order to demonstrate the effect of the explosive timescale existence to the system long-term evolution, the profiles of the tumor cell populations in the cases considered in Fig.~\ref{fig:RPSol+Tmscls} are depicted in Fig.~\ref{fig:RPSolComp}a.~In both cases, the initial tumor cell population $T(0)$ falls between the stable TFE and the unstable HTE.~Due to the existence of the explosive timescale at the initial period, the system is driven away from the stable TFE, as shown for both cases in the inset of Fig.~\ref{fig:RPSolComp}a.~At the tumor remission case, the explosive timescale disappears before $T$ reaches the unstable HTE threshold (P$_1$ in Fig.~\ref{fig:RPSolComp}) and as a result the system is driven towards the stable TFE, since only dissipative timescales exist.~On the other hand, in the tumor persistence case, the explosive timescale disappears after $T$ surpasses the unstable HTE threshold (P$_2$) and as a result, the system in this case is driven towards the stable HTE by the dissipative timescales.~Similar behavior applies in the tumor remission and persistence limiting cases considered in Fig.~\ref{fig:RPSolComp}b, in which $T(0)$ varies by only one cell.~As indicated by the inset of Fig.~\ref{fig:RPSolComp}b, the explosive timescale drives the system below or above the unstable HTE threshold, so that after losing its explosive nature (P$_3$ and P$_4$) the tumor is suppressed (stable TFE) or escapes immune surveillance (stable HTE).~It is thus, highlighted that the period in which an explosive timescale exists - denoted as \textit{explosive stage} from now on - is of particular interest, since it determines the long-term evolution of the system towards tumor suppression or remission.

%%%%%%%%%%%%%%%%%%%%%%%%%%%%%%%%%%%%%%%%%%%%%%%%%%%%%%%%%%%%%%%%%%%%%%%%%%%%%%%%
\section{The dynamics of the explosive stage}
\label{sec:CSPdyn}

It was established in Section~\ref{sec:exp} that the tumor-immune system dynamics is multi-scale in nature and exhibits an explosive stage, during which the long-term evolution of the system towards tumor persistence/remission is determined.

In this Section, the dynamics of the explosive stage is examined.~In order to generalize the results, a number of tumor persistence/remission cases were considered, among which detailed analysis is presented for 2 tumor persistence (TP and TP1) cases and 2 tumor remission (TR and TR1) ones; the related initial conditions are displayed in Table~\ref{tb:ICs}.~In particular, the TP case expresses the ``weakened" immune system of patient 9 in Ref.~\cite{de2006mixed}, which was adopted in Section~\ref{sec:exp}, while the TR case expresses a ``healthy" immune system being able to suppress a large tumor size.~In addition, the TP1 and TR1 cases considered in Fig.~\ref{fig:RPSolComp}b are examined as examples of limiting cases, since they only differ by 1 initial tumor cell.

\begin{table}[!h]
\centering
\begin{tabular}{l  c c c c  c c}
\toprule[1.5pt]
Case & T(0) & N(0) & L(0) & C(0) & long-term & $t_{exp}$ \\[2pt]
\midrule[1pt]
%TR & 10$^6$ & 10$^5$ & 10$^2$ & 6$\times$10$^{10}$ 	& TFE \\
TP & 10$^6$ & 10$^3$ & 10$^1$ & 6$\times$10$^{8}$ 	& HTE &	16.2\ d \\
TR &10$^7$ & 2$\times$10$^5$ & 10$^2$ & 4$\times$10$^{10}$ 	& TFE & ~2.3\ d	\\
TP1 & 319393 & 10$^3$ & 10$^1$ & 6$\times$10$^{8}$ 	& HTE &	37.3\ d\\
TR1 & 319392 & 10$^3$ & 10$^1$ & 6$\times$10$^{8}$ 	& TFE &	24.7\ d\\
\bottomrule[1.5pt]
\end{tabular}
\caption{The initial conditions of the 4 tumor persistence/remission cases considered and the corresponding duration of the explosive stage $t_{exp}$; parameter set of patient 9 \cite{de2006mixed}.}
\label{tb:ICs}
\end{table}

\begin{figure}[!b]
\centering
\subfigure[TP and TR cases]{\includegraphics[scale=0.23]{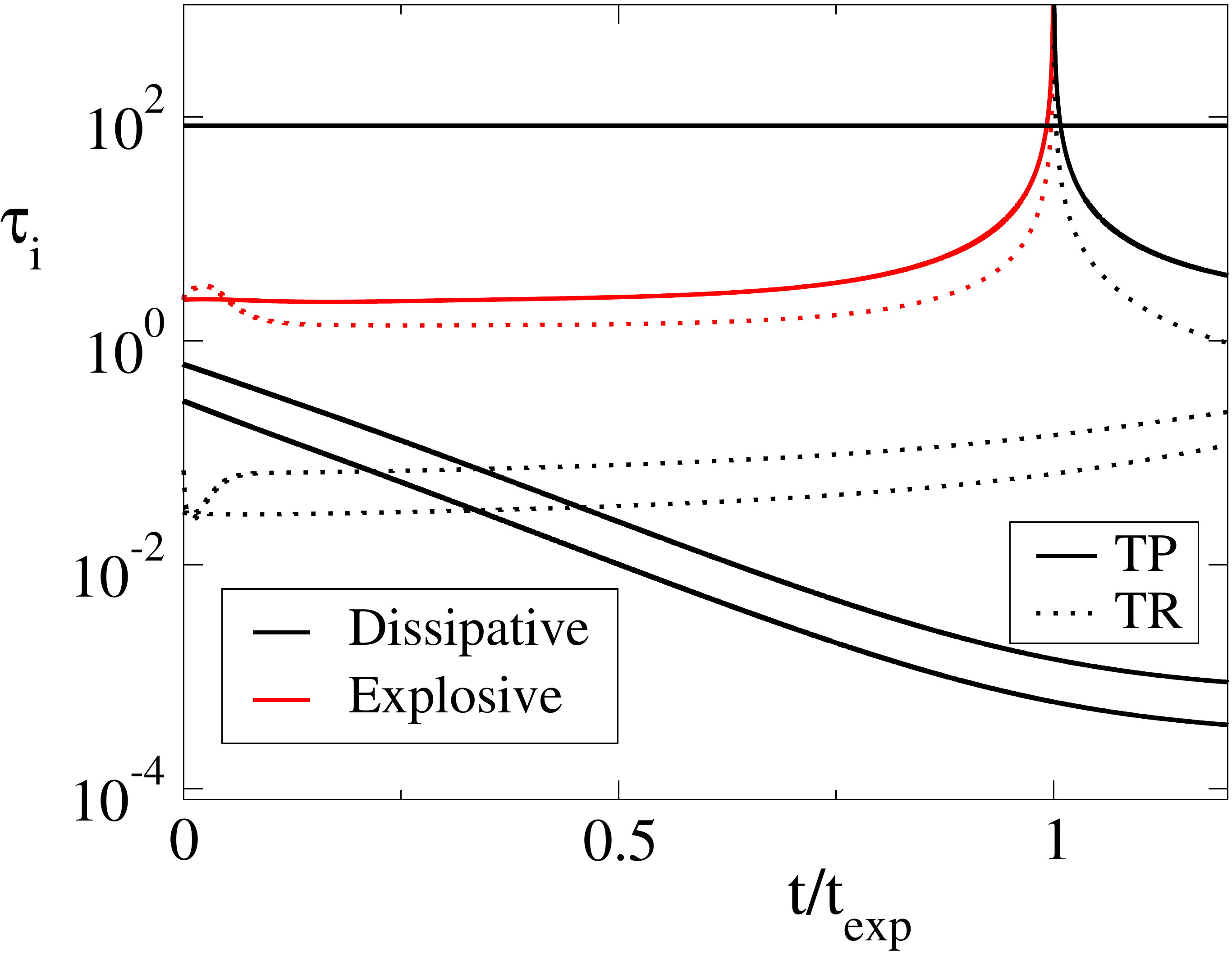}} \hspace{1.5cm} \subfigure[TP1 and TR1 cases]{\includegraphics[scale=0.23]{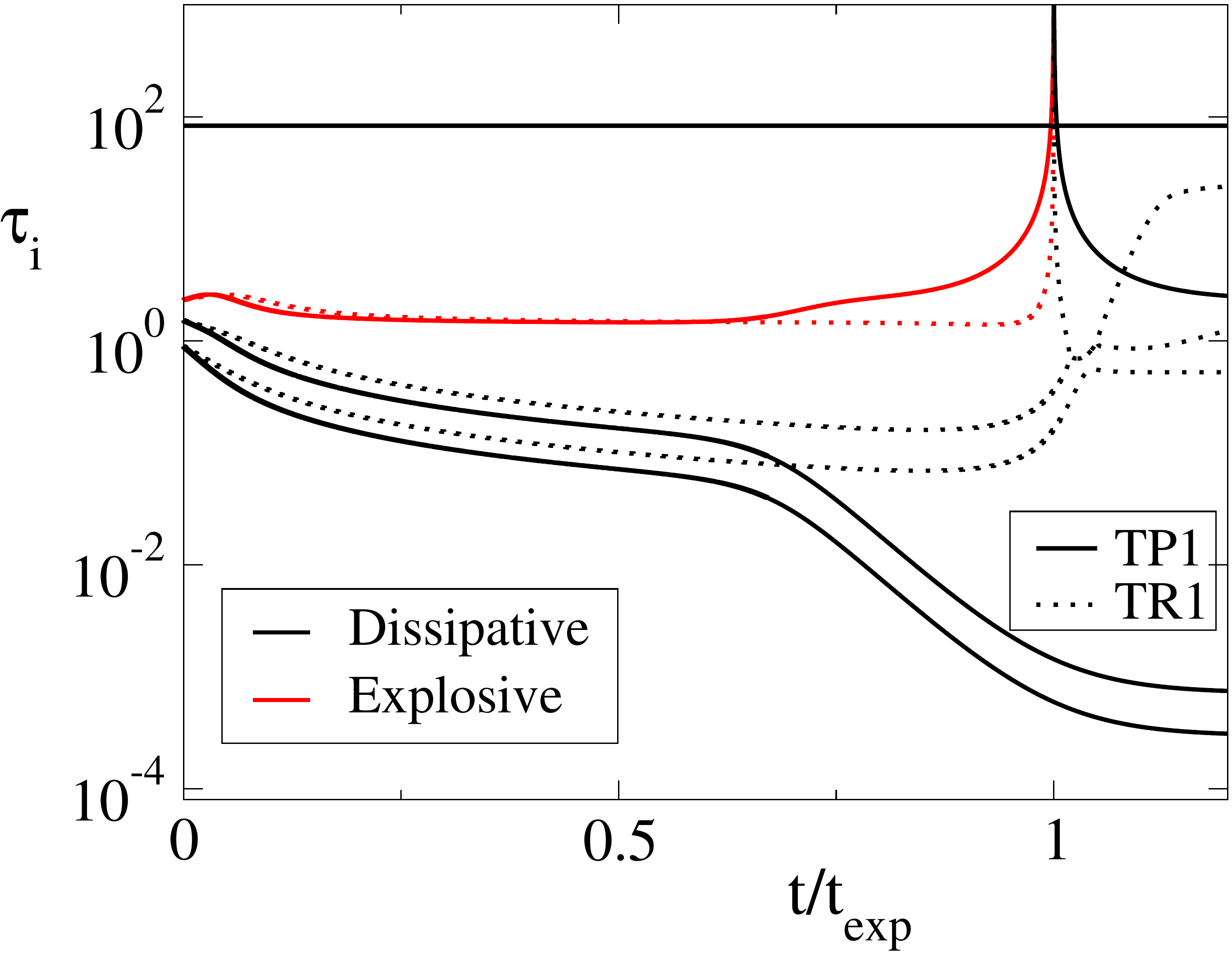}}
\caption{Timescales of the TP, TR, TP1 and TR1 cases normalized to $t_{exp}$.~The initial conditions and $t_{exp}$ of each case is displayed in Table~\ref{tb:ICs}.}
\label{fig:Tmscls_cases}
\end{figure}

As discussed in Section~\ref{sec:exp}, the system exhibits 4 timescales, the evolution of which is displayed in Fig.~\ref{fig:Tmscls_cases} normalized to the corresponding $t_{exp}$ for all cases considered.~The two fastest timescales $\tau_{1,2}$ are dissipative in nature, tending to drive the system towards a stable equilibrium, either HTE or TFE.~The 3rd timescale $\tau_3$ is explosive at the explosive stage and becomes dissipative only after it, so that its action initially drives the system away from the stable equilibrium and towards it (the same equilibrium or differennt one) after the explosive stage.~In all cases considered in Fig.~\ref{fig:Tmscls_cases}, the gap between $\tau_2$ and $\tau_3$ becomes significant after a rapid transient period, indicating the strong multi-scale nature of the system.~At the tumor persistence TP and TP1 cases, this gap is wider, indicating a more pronounced multi-scale character as denoted in Section~\ref{sec:exp}.~In addition, the gap between $\tau_2$ and $\tau_3$ indicates that during the explosive stage $\tau_{1,2}$ become quickly exhausted after a rapid transient period, since they are both dissipative in nature and much faster than the explosive timescale $\tau_3$.

Theerfore, during the explosive stage the system is multi-scale in nature (fast/slow timescale seperation) and exhibits $M=2$ exhausted and dissipative timescales, which are much faster than the slow and explosive $\tau_3$.~Under these conditions, according to GSPT \citep{fenichel1979geometric,kaper1999systems}, the system is confined to evolve along a $M=2$-dim. SIM under the action of the fast dissipative $\tau_{1,2}$.~The slow evolution of the system along the SIM is then governed by a slow system, which is characterized by action of the slow timescales, mainly by the fastest of the slow timescales $\tau_3$ being explosive during the explosive stage~\citep{verhulst2006nonlinear,kuehn2015multiple}.~In order to identify the components of the model in Eq.~\eqref{eq:VecF} that form the constraints giving rise to the 2-dim. SIM and govern the slow evolution of the system, the CSP methodology and its algorithmic diagnostics tools \citep{lam1994,lam1989} were employed; a brief presentation of CSP and the related diagnostics is provided in \ref{app:CSP}.

%%%%%%%%%%%%%%
\subsection{The underlying mechanism of tumor progression during the explosive stage}
\label{subsec:fast}

First, the underlying mechanisms of tumor progression are revealed, by identifying the constraints along which the system is bound to evolve during the explosive stage.~According to CSP in \ref{app:CSP}, the system in Eq.~\eqref{eq:VecF} can be decomposed in 4 CSP modes, the first $M=2$ of which
become exhausted after a rapid transient period at the explosive stage, thus imposing the emergence of $M=2$ constraints.~Employing the CSP diagnostic tools for each exhausted CSP mode enables the algorithmic identification of (a) the processes contributing to the formation of each constraint via the \textit{Amplitude Participation Index} (API) in Eq.~\eqref{eq:gov5}, (b) the processes driving the emergence of each constraint via the \textit{Timescale Participation Index} (TPI) in Eq.~\eqref{eq:gov6}, and (c) the cell populations mostly related to each constraint via the \textit{CSP Pointer} (Po) in Eq.~\eqref{eq:gov7}.~Details on the formulation of each CSP tool is provided in \ref{appSub:tools}.~The related API, TPI and Po identifications are displayed in Table~\ref{tb:M1diags} for each exhausted CSP mode at the four TP, TR, TP1 and TR1 cases considered.~Being interested to the explosive stage, the corresponding identifications are reported in Table~\ref{tb:M1diags} at representative points during the explosive stage at $t/t_{exp}=0.2, 0.5, 0.8$.~Note that during the initial transient period, the constraints have not been established yet, as indicated by the CSP identifications at the beginning of the process $t/t_{exp}=0$ in Table~\ref{tb:M1diags}.

During the explosive stage, according to the API values in Table~\ref{tb:M1diags}, the cancellation of process 2 and process 13 indicate the formation of the 1st constraint, so that
\begin{equation}
R^2  \approx R^{13}   \Rightarrow   e C   \approx  p N T 
\label{eq:Con1}
\end{equation}
which expresses the equilibration between the growth rate of NK cells $R^2$ and the inactivation rate of NK cells $R^{13}$.~According to TPI, the formation of the 1st constraint in Eq.~\eqref{eq:Con1} is mainly driven by the fast process 13 and is exclusively related to the population of NK cells, $N$, as indicated by the Po values in Table~\ref{tb:M1diags}.

In addition, the formation of the 2nd constraint is indicated by the cancellation of processes 12 with 14 and 5, according to the API values in Table~\ref{tb:M1diags}, so that
\begin{equation}
R^{12} \approx  R^{14}+R^5   \Rightarrow   r_2 C T   \approx  q L T + m L
\label{eq:Con2}
\end{equation}
which expresses the equilibration between the recruitment rate of the CD8$^+$ T cells $R^{12}$ and the inactivation rate of CD8$^+$ T cells $R^{14}$, accompanied by a small contribution by the death rate of CD8$^+$ T cells $R^5$, which is more pronounced at the TP1 and TR1 cases.~According to TPI, the formation of the 2nd constraint in Eq.~\eqref{eq:Con2} is mainly driven by the fast process 14; much lesser contributions are provided by processes 5 and 12 mostly at the TP1 and TR1 cases.~Note however, that processes 14 and 5 contribute towards the formation of the 2nd constraint (negative TPI values in Table~\ref{tb:M1diags}), while process 12 opposes to its formation (positive TPI values).~This is because, as indicated by Po this constraint is exclusively related to the population of CD8$^+$ T cells, $L$, and process 12 (recruitment) tends to produce CD8$^+$ T cells, while processes 14 and 5 (inactivation and death) tend to consume them.

\begin{table}[!h]
\scriptsize
\centering
\begin{tabular}{c c | lr  lr  lr  lr } 
& \multicolumn{9}{c}{~~~~~~~~~1st CSP mode} \\
\toprule[1.5pt]
& $t/t_{exp}$	& \multicolumn{2}{c}{0} & \multicolumn{2}{c}{0.2} & \multicolumn{2}{c}{0.5} & \multicolumn{2}{c}{0.8}  \\[2pt] \midrule[1pt]
\parbox[c]{1mm}{\multirow{8}{*}{\rotatebox[origin=c]{90}{\textbf{TP case} }}} & \multirow{4}{*}{API} 	&	13	&	0.85	&	\tikzmark{top left 1}13	&	0.49	&	13	&	0.50	&	13	&	0.50\tikzmark{top right 1}	\\
& 				&	1	&	-0.09	&	\tikzmark{bottom left 1}2	&	-0.48	&	2	&	-0.50	&	2	&	-0.50\tikzmark{bottom right 1}	\\
&				&	2	&	-0.03	&	1	&	-0.02	&		&		&		&			\\
& & & & & & & & & \\[2pt] \cmidrule{2-10}
& \multirow{2}{*}{TPI} 	&	13	&	-0.99	&	13	&	-1.00	&	13	&	-1.00	&	13	&	-1.00		\\
& & & & & & & & & \\[2pt] \cmidrule{2-10}
& Po	&	$N$	&	1.00	&	$N$	&	1.00	&	$N$	&	1.00	&	$N$	&	1.00	\\[2pt] \midrule[1pt]
\parbox[c]{1mm}{\multirow{10}{*}{\rotatebox[origin=c]{90}{\textbf{TR case} }}} & \multirow{4}{*}{API}	&	13	&	0.81	&	\tikzmark{top left 4}2	&	-0.49	&	2	&	-0.48	&	2	&	-0.47\tikzmark{top right 4}		\\
&				&	12	&	-0.14	&	\tikzmark{bottom left 4}13	&	0.48	&	13	&	0.48	&	13	&	0.46\tikzmark{bottom right 4}		\\
&				&	14	&	0.02	&		&		&		&		&	8	&	0.02		\\
&				&	15	&	0.01	&		&		&		&		&		&			\\[2pt] \cmidrule{2-10}
& \multirow{3}{*}{TPI}	&	13	&	-0.43	&	13	&	-1.00	&	13	&	-1.00	&	13	&	-1.00		\\
&				&	15	&	-0.33	&		&		&		&		&		&				\\
&				&	14	&	-0.23	&		&		&		&		&		&				\\[2pt] \cmidrule{2-10}
& \multirow{2}{*}{Po}	&	$L$	&	0.56	&	$N$	&	1.00	&	$N$	&	1.00	&	$N$	&	1.00		\\
&				&	$N$	&	0.44	&		&		&		&		&		&			\\[2pt] \midrule[1pt]
\parbox[c]{1mm}{\multirow{10}{*}{\rotatebox[origin=c]{90}{\textbf{TP1 case} }}} & \multirow{5}{*}{API}	&	13	&	0.64	&	\tikzmark{top left 2}13	&	0.44		&	13	&	0.47	&	13	&	0.50\tikzmark{top right 2}	\\
&				&	1	&	-0.18	&	\tikzmark{bottom left 2}2	&	-0.44		&	2	&	-0.47	&	2	&	-0.50\tikzmark{bottom right 2}	\\
&				&	3	&	0.08	&	12	&	-0.03	&	12	&	-0.02	&		&			\\
&				&	2	&	-0.07	&	14	&	0.03	&	14	&	0.02	&		&			\\
&				&		&		&	1	&	-0.03	&		&		&		&			\\[2pt] \cmidrule{2-10}
& \multirow{3}{*}{TPI}				&	13	&	-0.95	&	13	&	-0.99	&	13	&	-1.00	&	13	&	-1.00		\\
& & & & & & & & & \\
& & & & & & & & & \\[2pt] \cmidrule{2-10}
& Po				&	$N$	&	1.00	&	$N$	&	1.00	&	$N$	&	1.00	&	$N$	&	1.00\\[2pt] \midrule[1pt]
\parbox[c]{1mm}{\multirow{10}{*}{\rotatebox[origin=c]{90}{\textbf{TR1 case} }}} & \multirow{5}{*}{API}	&	13	&	0.64	&	\tikzmark{top left 3}2	&	-0.43	&	13	&	0.46	&	13	&	0.47\tikzmark{top right 3}		\\
&				&	1	&	-0.18	&	\tikzmark{bottom left 3}13	&	0.42	&	2	&	-0.46	&	2	&	-0.47\tikzmark{bottom right 3}	\\
&				&	3	&	0.08	&	12	&	-0.04	&	12	&	-0.02	&	12	&	-0.02 \\
&				&	2	&	-0.07	&	1	&	-0.03	&	14	&	0.02	&	14	&	0.02	\\
&				&		&		&	14	&	0.03	&		&		&		&		\\[2pt] \cmidrule{2-10}
& \multirow{3}{*}{TPI}	&	13	&	-0.95	&	13	&	-0.99	&	13	&	-0.99	&	13	&	-1.00 	\\
& & & & & & & & & \\
& & & & & & & & & \\[2pt] \cmidrule{2-10}
& Po				&	$N$	&	1.00	&	$N$	&	1.00	&	$N$	&	1.00	&	$N$	&	1.00		\\[2pt] 
\bottomrule[1.5pt]
\end{tabular}
\begin{tabular}{lr  lr  lr  lr  } 
\multicolumn{8}{c}{2nd CSP mode} \\
\toprule[1.5pt]
\multicolumn{2}{c}{0} & \multicolumn{2}{c}{0.2} & \multicolumn{2}{c}{0.5} & \multicolumn{2}{c}{0.8}  \\[2pt] \midrule[1pt]
12	&	0.50	&	\tikzmark{top left 5}12	&	0.50	&	12	&	0.50	&	12	&	0.50\tikzmark{top right 5} 	\\
3	&	-0.38	&	\tikzmark{bottom left 5}14	&	-0.46	&	14	&	-0.50	&	14	& -0.50\tikzmark{bottom right 5}	\\
1	&	-0.10	&	3	&	-0.02	&		&		&		&				\\
14	&	-0.01	&	5	&	-0.02	&		&		&		&				\\[2pt] \midrule
14	&	-0.87	&	14	&	-0.96	&	14	&	-1.00	&	14	&	-1.00		\\
5	&	-0.13	&	5	&	-0.04	&		&		&		&				\\[2pt] \midrule
$L$	&	1.00	&	$L$	&	1.00	&	$L$	&	1.00	&	$L$	&	1.00		\\[2pt] \midrule[1pt]
13	&	0.84	&	\tikzmark{top left 8}12	&	0.49	&	12	&	0.49	&	12	&	0.48\tikzmark{top right 8}		\\
12	&	0.12	&	\tikzmark{bottom left 8}14	&	-0.48	&	14	&	-0.48	&	14	&	-0.46\tikzmark{bottom right 8}	\\
	&		&		&		&		&		&	2	&	-0.02		\\
	&		&		&		&		&		&	13	&	0.02			\\[2pt] \midrule
13	&	-0.55	&	14	&	-0.92	&	14	&	-0.90	&	14	&	-0.85\\
15	&	-0.25	&	12	&	0.05	&	12	&	0.07	&	12	&	0.09		\\
14	&	-0.18	&		&		&		&		&		&		\\[2pt] \midrule
$N$	&	0.56	&	$L$	&	1.00	&	$L$	&	1.00	&	$L$	&	1.01		\\
$L$	&	0.44	&		&		&		&		&		&			\\[2pt] \midrule[1pt]
3	&	-0.57	&	\tikzmark{top left 6}12	&	0.49	&	12	&	0.49	&	12	&	0.50\tikzmark{top right 6}	\\
12	&	0.30	&	14	&	-0.43	&	14	&	-0.47	&	14	&	-0.50		\\
1	&	-0.12	&	\tikzmark{bottom left 6}5	&	-0.03	&	5	&	-0.02	&		&		\tikzmark{bottom right 6} \\
& & & & & & &  \\
& & & & & & &  \\[2pt] \midrule
14	&	-0.69	&	14	&	-0.80	&	14	&	-0.88	&	14	&	-1.00		\\
5	&	-0.31	&	12	&	0.13	&	12	&	0.08	&		&				\\
	&		&	5	&	-0.05	&	5	&	-0.03	&		&				\\[2pt] \midrule
$L$	&	1.00	&	$L$	&	1.02 	&	$L$	&	1.00	&	$L$	&	1.00		\\[2pt] \midrule[1pt]
3	&	-0.57	&	\tikzmark{top left 7}12	&	0.48	&	12	&	0.49	&	12	&	0.49\tikzmark{top right 7}	\\
12	&	0.30	&	14	&	-0.40	&	14	&	-0.46	&	14	&	-0.47		\\
1	&	-0.12	&	\tikzmark{bottom left 7}5	&	-0.04	&	5	&	-0.02	&	5	&	-0.02\tikzmark{bottom right 7}		\\
&		&	3	&	-0.04	&		&		&		&		\\
& & & & & & &  \\[2pt] \midrule
14	&	-0.69	&	14	&	-0.76	&	14	&	-0.85	&	14	&	-0.87		\\
5	&	-0.31	&	12	&	0.14	&	12	&	0.10	&	12	&	0.08		\\
	&		&	5	&	-0.07	&	5	&	-0.04	&	5	&	-0.03		\\[2pt] \midrule
$L$	&	1.00	&	$L$	&	1.03	&	$L$	&	1.01	&	$L$	&	1.00		\\[2pt]
 \bottomrule[1.5pt]
\end{tabular}
\DrawBox[ultra thick, white, fill=black, fill opacity=0.2]{top left 1}{bottom right 1}
\DrawBox[ultra thick, white, fill=black, fill opacity=0.2]{top left 2}{bottom right 2}
\DrawBox[ultra thick, white, fill=black, fill opacity=0.2]{top left 3}{bottom right 3}
\DrawBox[ultra thick, white, fill=black, fill opacity=0.2]{top left 4}{bottom right 4}
\DrawBox[ultra thick, white, fill=black, fill opacity=0.2]{top left 5}{bottom right 5}
\DrawBox[ultra thick, white, fill=black, fill opacity=0.2]{top left 6}{bottom right 6}
\DrawBox[ultra thick, white, fill=black, fill opacity=0.2]{top left 7}{bottom right 7}
\DrawBox[ultra thick, white, fill=black, fill opacity=0.2]{top left 8}{bottom right 8}
\caption{The API, TPI and Po of the 1st and 2nd CSP modes for the TP, TR, TP1 and TR1 cases during the explosive stage.~The major contributions are presented at the beginning of the process $t/t_{exp}=0$ and at $t/t_{exp}=0.2, 0.5, 0.8$ of the explosive stage (skewed boxes); sums of the absolute API, TPI and Po values up to more than 95\% (0.95).}
\label{tb:M1diags}
\end{table}

The processes forming the constraints in Eqs.~\eqref{eq:Con1} and \eqref{eq:Con2}, the ones driving their formation and the cell populations related to them generalize well for all the cases considered.~Despite the lesser but still major contribution of processes 2 and 13 to the 1st CSP mode at the TP1 and TR1 cases and the minor contribution of process 5 to the 2nd CSP mode at the TP and TR cases, the CSP identifications indicate that the system is confined to evolve along these constraints throughout the explosive stage, after a small initial transient period; compare the identifications at skewed boxes with the ones at $t/t_{exp}=0$ in Table~\ref{tb:M1diags}.

\begin{figure}[!h]
\centering
\subfigure[TP case profiles]{\includegraphics[scale=0.18]{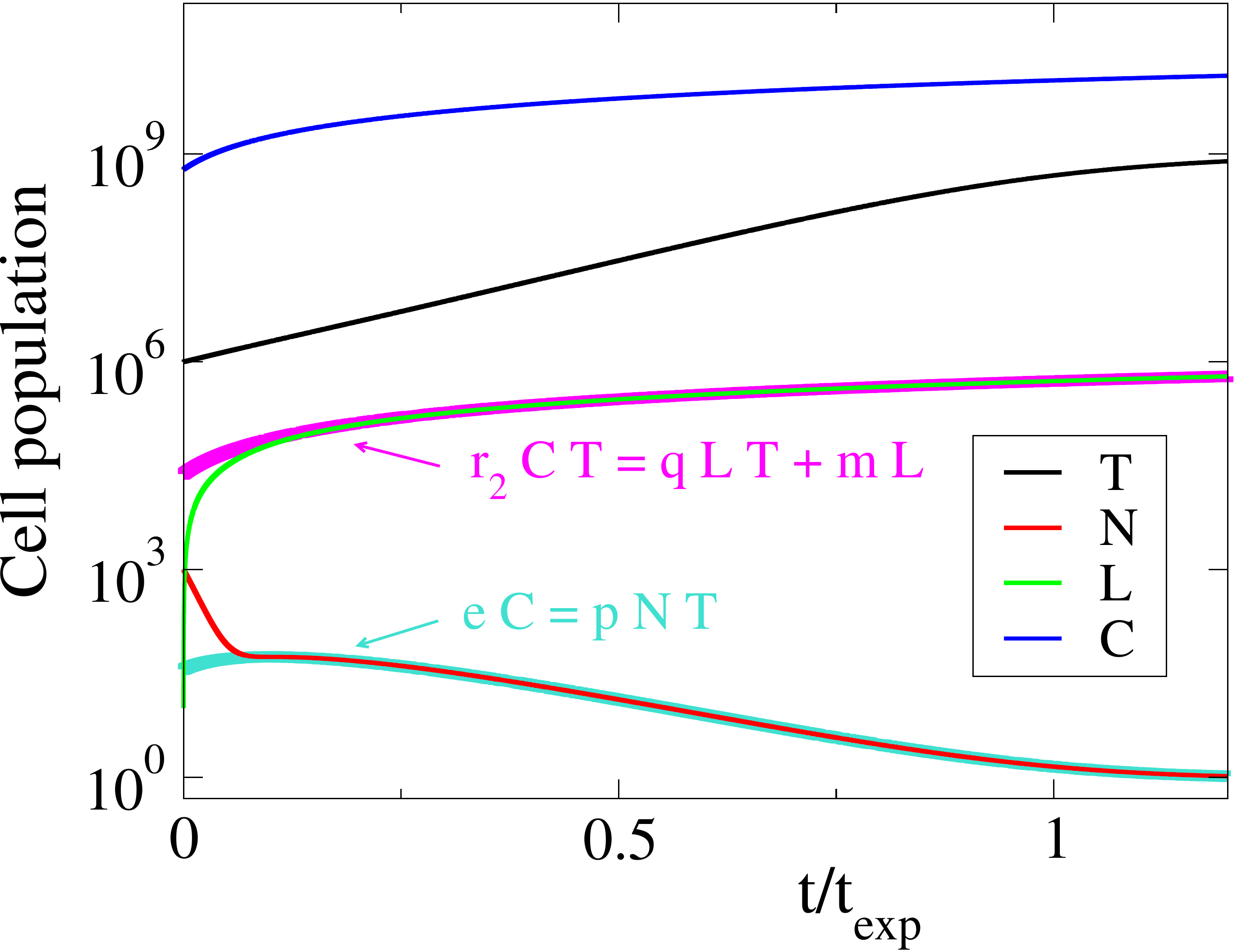}} \hspace{0.05cm} \subfigure[$RE_N=(N-\hat{N})/N$]{\includegraphics[scale=0.18]{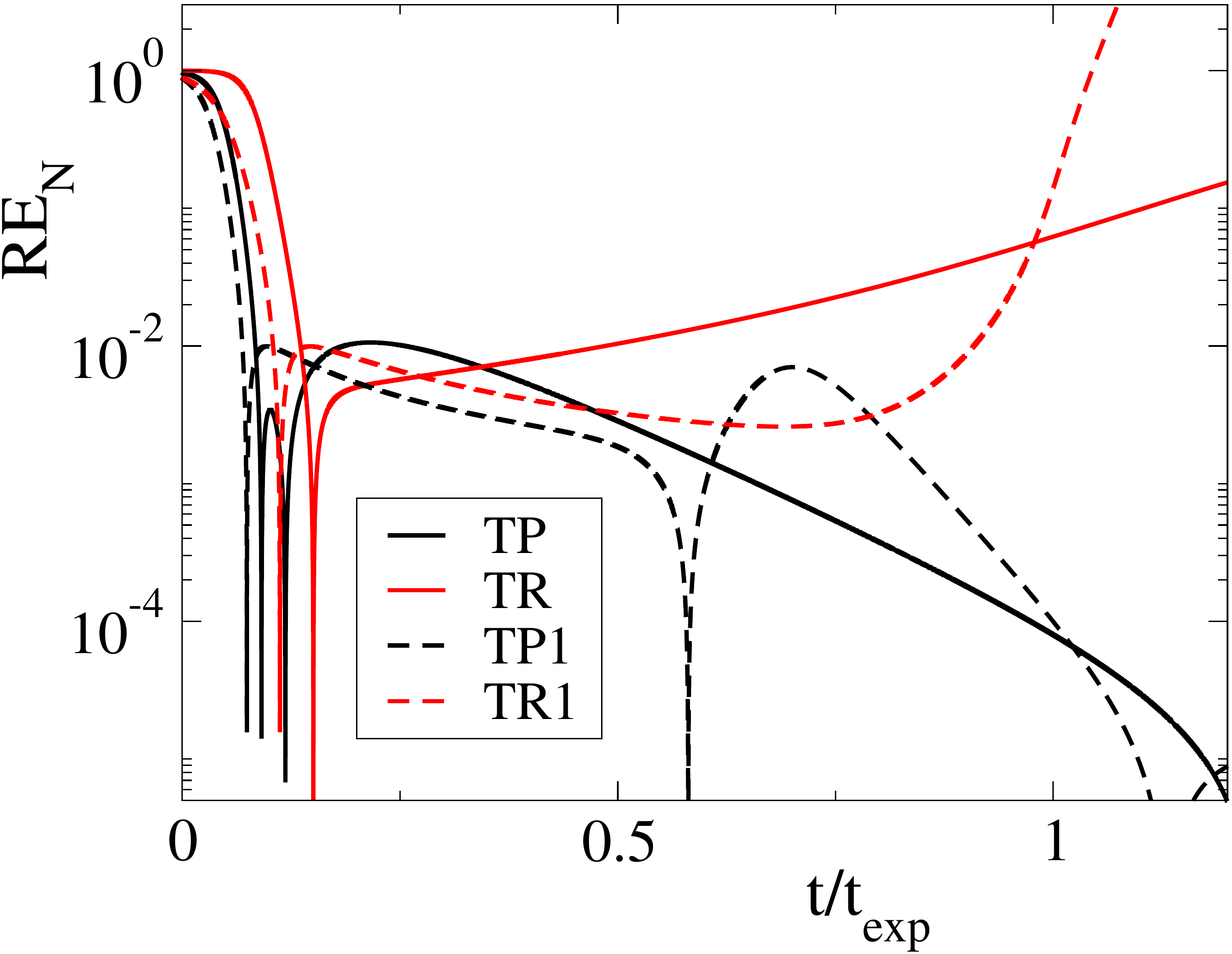}} \hspace{0.05cm} \subfigure[$RE_L=(L-\hat{L})/L$]{\includegraphics[scale=0.18]{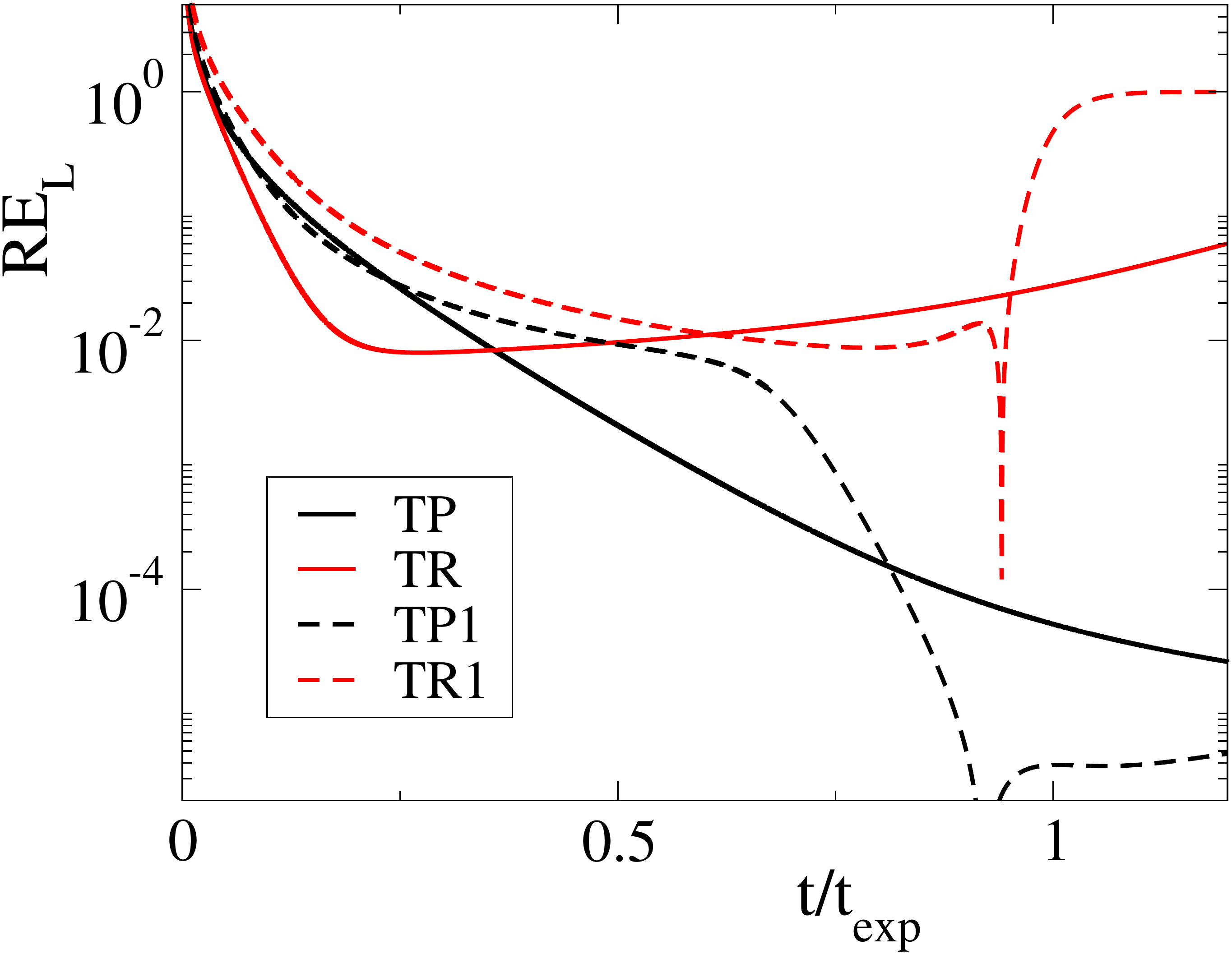}}
\caption{The profiles of $N$ and $L$ on the basis of the system in Eq.~\eqref{eq:VecF} compared to their calculation via the algebraic relations of the constraints in Eqs.~\eqref{eq:Con1} and \eqref{eq:Con2}, respectively; cyan and magenta thick lines on the left panel.} 
\label{fig:TP_Sol+SIM_err9}
\end{figure}

In order to validate the CSP identifications, the profiles of $N$ and $L$ obtained on the basis of the full model in Eq.~\eqref{eq:VecF} were compared to the expressions $\hat{N}$ and $\hat{L}$ provided by the algebraic equations of the constraints in Eqs.~\eqref{eq:Con1} and \eqref{eq:Con2}.~In particular, as shown in Fig.~\ref{fig:TP_Sol+SIM_err9}a for the TP case, after the initial transient period, $N$ and $L$ are perfectly aligned to the calculated by Eqs.~\eqref{eq:Con1} and \eqref{eq:Con2} expressions of $\hat{N}=eC/pT$ (cyan line) and $\hat{L}=r_2CT/(qT+m)$ (magenta line), respectively.~According to the relative error $RE_N=(N-\hat{N})/N$ and $RE_L=(L-\hat{L})/L$ that is depicted for all the cases considered for patient 9 parameter set in Fig.~\ref{fig:TP_Sol+SIM_err9}b,c, the full model solution of $N$ and $L$ is approximated with great accuracy by $\hat{N}$ and $\hat{L}$, during the explosive stage.~In addition, the validity of the constraints in Eqs.~\eqref{eq:Con1} and \eqref{eq:Con2} is preserved when considering different parameter values as shown in Fig.~S1, for which the parameter set of patient 10 in Ref.~\cite{de2006mixed} was adopted at 4 indicative tumor progression cases.~Finally note that the constraint-derived solutions of $\hat{N}$ and $\hat{L}$ approximate accurately the profiles of the full model $N$ and $L$ at all the TP cases considered for both patient 9 and 10 parameter sets, both during the explosive stage and after it.~This feature indicates that the constraints in Eqs.~\eqref{eq:Con1} and \eqref{eq:Con2} are preserved after the explosive stage only at the TP cases.

%%%%%%%%%%%%%%
\subsection{The governing dynamics during the explosive stage}
\label{subsec:slow}

Having obtained the insight on the constraints along which the system is confined to evolve during the explosive stage, next the focus is turned on the slow dynamics of the system along these constraints.~The slow evolution of the system is governed by a slow system,  which consists of the remaining $N-M=2$ slow CSP modes and is governed by the \textit{characteristic} - the fastest active - timescale, that is $\tau_3=\tau_{exp}$ during the explosive stage.~Employing the CSP diagnostic tools for the slow CSP modes enables the algorithmic identification of (a) the processes driving the dynamics of the slow system via the TPI of the characteristic 3rd CSP mode in Eq.~\eqref{eq:gov6} and (b) the cell populations mostly mostly related to the slow system dynamics via the Po of the characteristic 3rd CSP mode in Eq.~\eqref{eq:gov7}.~The related TPI and Po identifications are displayed in Table~\ref{tb:M3diags} for the characteristic CSP mode at the four TP, TR, TP1 and TR1 cases considered.~Similarly to Table~\ref{tb:M1diags}, the corresponding identifications are reported in Table~\ref{tb:M3diags} at representative points during the explosive stage at $t/t_{exp}=0.2, 0.5, 0.8$.~Note that the identifications of the 3rd CSP mode at $t/t_{exp}=0$ in Table~\ref{tb:M3diags} do not provide insights on the slow dynamics of the system, since the constraints have not been established yet during the initial transient period.

\begin{table}[!t]
\scriptsize
\centering
\begin{tabular}{c c | lr  lr  lr  lr} \toprule[1.5pt]
& $t/t_{exp}$ 	& \multicolumn{2}{c}{0} & \multicolumn{2}{c}{0.2} & \multicolumn{2}{c}{0.5} & \multicolumn{2}{c}{0.8}  \\[2pt] \midrule[1pt]
\parbox[c]{3mm}{\multirow{6}{*}{\rotatebox[origin=c]{90}{\textbf{TP case} }}}  & \multirow{4}{*}{TPI}&	1	&	1.00	&	1	&	0.80	&	1	&	0.98	&	1	&	1.00	\\
& &	&		&	12	&	-0.08	&	12	&	-0.01	&		&		\\
& &	&		&	14	&	0.08	&	14	&	0.01	&		&		\\
& &	&		&	8	&	0.04	&		&		&		&		\\[2pt] \cmidrule{2-10}
& Po	& $T$	 &	1.00	&	$T$	&	1.00	&	$T$	&	1.00	&	$T$	&	1.00		\\[2pt] \midrule[1pt]
\parbox[c]{3mm}{\multirow{6}{*}{\rotatebox[origin=c]{90}{\textbf{TR case} }}}  & \multirow{4}{*}{TPI}	&	1	&	0.68	&	12	&	-0.35	&	12	&	-0.36	&	12	&	-0.40		\\
&	&	12	&	-0.12	&	14	&	0.35	&	14	&	0.36	&	14	&	0.40		\\
&	&	8	&	-0.09	&	1	&	0.17	&	1	&	0.16	&	1	&	0.15		\\
&	&	14	&	0.06	&	8	&	0.12	&	8	&	0.11	&	8	&	0.04		\\[2pt] \cmidrule{2-10}
&	Po &	$T$	&	1.00	&	$T$	&	1.00	&	$T$	&	1.00	&	$T$	&	1.01		\\[2pt] \midrule[1pt]
\parbox[c]{3mm}{\multirow{6}{*}{\rotatebox[origin=c]{90}{\textbf{TP1 case} }}}  & \multirow{4}{*}{TPI} &	1	&	1.00	&	12	&	-0.31	&	12	&	-0.32	&	1	&	0.88		\\
&	&		&		&	14	&	0.30	&	14	&	0.32	&	14	&	0.05		\\
&	&		&		&	1	&	0.27	&	1	&	0.23	&	12	&	-0.05		\\
&	&		&		&	8	&	0.13	&	8	&	0.14	&	8	&	0.02	\\[2pt] \cmidrule{2-10}
& Po &	$T$	&	1.00	&	$T$	&	1.02	&	$T$	&	1.00	&	$T$	&	1.00		\\[2pt] \midrule[1pt]
\parbox[c]{3mm}{\multirow{6}{*}{\rotatebox[origin=c]{90}{\textbf{TR1 case} }}}  & \multirow{4}{*}{TPI}	&	1	&	1.00	&	1	&	0.30	&	12	&	-0.31	&	12	&	-0.32		\\
&	&		&		&	12	&	-0.30	&	14	&	0.31	&	14	&	0.32		\\
&	&		&		&	14	&	0.29	&	1	&	0.24	&	1	&	0.22	\\
&	&		&		&	8	&	0.11	&	8	&	0.13	&	8	&	0.14		\\[2pt] \cmidrule{2-10}
&	Po	&	$T$	&	1.00	&	$T$	&	1.03	&	$T$	&	1.01	&	$T$	&	1.01		\\[2pt]  \bottomrule[1.5pt]
\end{tabular}
%\DrawBox[ultra thick, white, fill=black, fill opacity=0.2]{top left 1}{bottom right 1}
\caption{The TPI and Po of the characteristic, during the explosive stage, 3rd mode for the TP, TR, TP1 and TR1 cases.~The major contributions are presented at the beginning of the process $t/t_{exp}=0$ and at $t/t_{exp}=0.2, 0.5, 0.8$ of the explosive stage; sums of the absolute API, TPI and Po values up to more than 95\% (0.95).}
\label{tb:M3diags}
\end{table}

According to Table~\ref{tb:M3diags}, the processes driving the slow system during the explosive stage are:
\begin{enumerate}[(i)]
\item primarily the tumor growth (process 1), since it is consistently identified by the TPI to majorly promote the explosive character of $\tau_{exp}$ (positive TPI).~The role of tumor growth process is more pronounced at the tumor persistence TP and TP1 cases, and its contribution on driving the slow system becomes more significant towards the end of the explosive stage.~On the other hand, at the tumor remission TR and TR1 cases, the tumor growth process provides lesser contributions in driving the slow system, having however still significant importance.
\item the recruitment and inactivation of CD8$^+$ T cells (processes 12 and 14); the former opposing to the explosive character of $\tau_{exp}$ and the latter promoting it (negative and positive TPIs).~Note that these processes were identified to equilibrate for the formation of the 3nd constraint in Eq.~\eqref{eq:Con2}, explaining their contrasting contributions in driving the slow system.~The contribution of processes 12 and 14 are more pronounced at the tumor remission TR and TR1 cases, while they become of lesser importance towards the end of the explosive stage at the tumor persistence TP and TP1 cases.
\item to a lesser degree the fractional cell kill rate of tumor cells by CD8$^+$ T cells (process 8).~This process, despite tending to decrease the tumor cell population, promotes the explosive character of $\tau_{exp}$ and becomes less significant towards the end of the explosive stage; compare the TPI of process 8 in Table~\ref{tb:M3diags} along the explosive stage.
\end{enumerate}

In addition, according to Po in Table~\ref{tb:M3diags} for all the cases considered, the tumor cell population is mostly related to the slow system dynamics during the explosive stage; a result that is in agreement with the driving processes 1 and 8 that mainly produce and consume tumor cells, respectively.~This result indicates that the slow dynamics primarily affects the slow evolution of tumor cell population.

\subsubsection{Introducing perturbations on the parameter values}
\label{subsubsec:demo}

\begin{figure}[!t]
\centering
\subfigure[80\% and 120\% perturbations on $a$]{\includegraphics[scale=0.22]{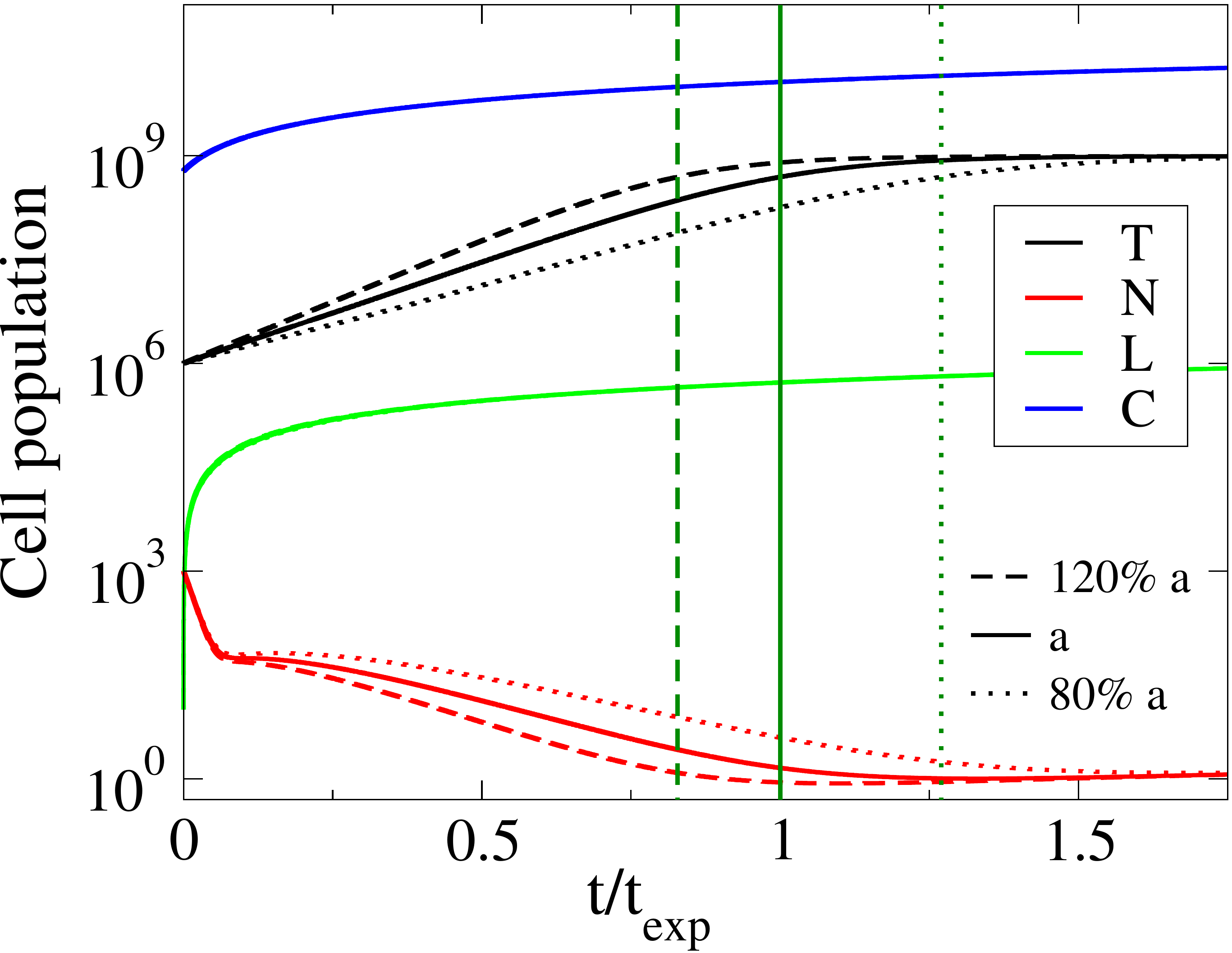}} \hspace{1.5cm} \subfigure[40\% perturbation on $c$ and $e$]{\includegraphics[scale=0.22]{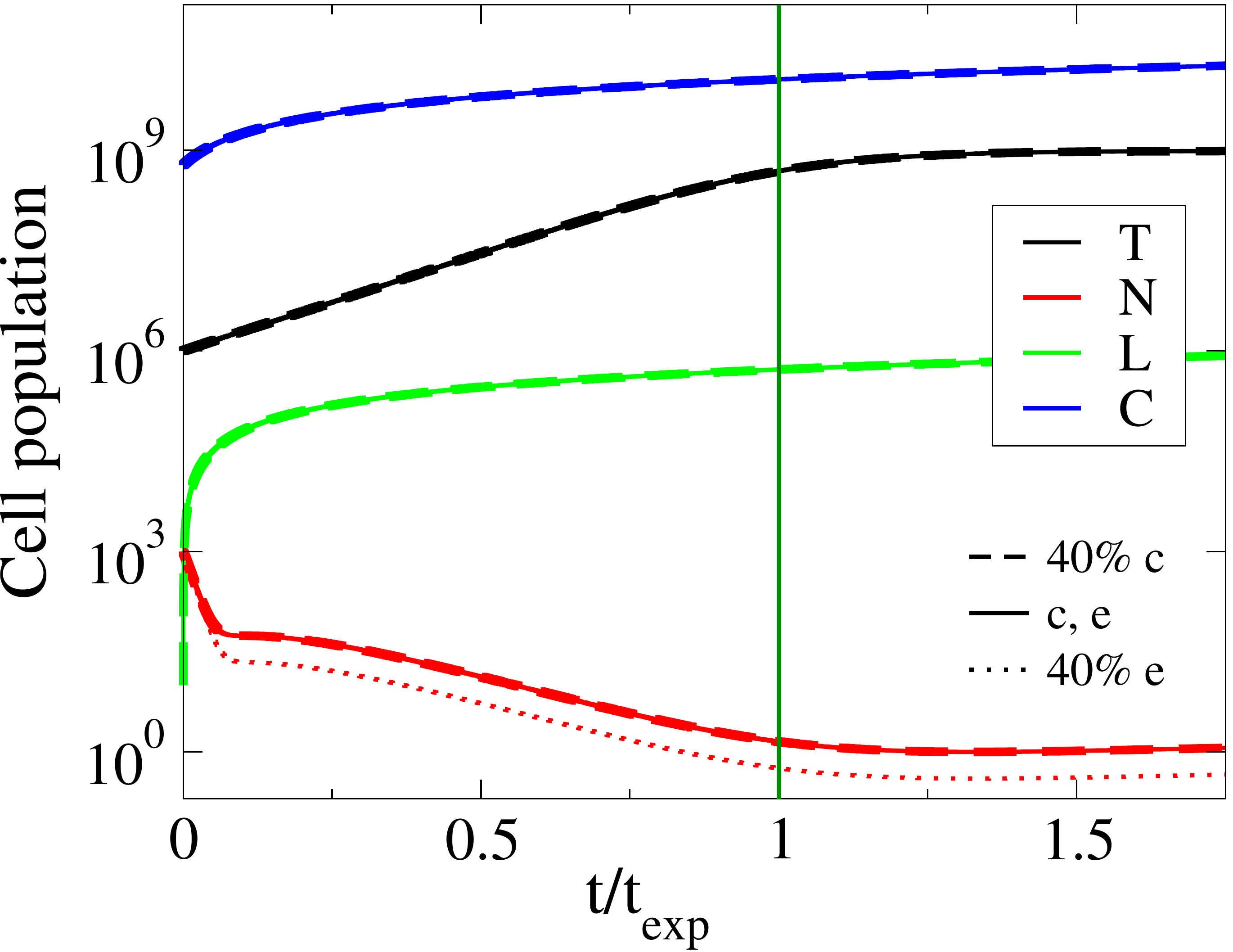}}
\caption{Cell population responses in perturbations on $a$, $c$ and $e$ in comparison to the nominal tumor persistence case.~The dark green vertical lines indicate the end of the explosive stage for each (un)perturbed case.}
\label{fig:1_abP}
\end{figure}

In order to validate the above conclusions related to the governing slow dynamics, the response of the system in various perturbations is examined.~In particular, the parameters $a$, $c$ and $e$, the rate constants of processes 1, 7 and 2 respectively, are perturbed in order to demonstrate the response of the governing slow dynamics; the response of the system in the above perturbations is depicted in Fig.~\ref{fig:1_abP}.

Starting with the process majorly driving the slow system, a potential reduction of $a$ will result in reducing the rate of process 1, that is the process majorly generating $\tau_{exp}$.~As a result $\tau_{exp}$ will become slower, so that the explosive stage will last longer.~In addition, the reduction of process 1 will further result in reducing the population of tumor cells, since this is the cell population mostly mostly related to the slow system dynamics during the explosive stage, according to Po in Table~\ref{tb:M3diags}.~In turn, the reduction of $T$ will result to an increase of $N$ and a minor decrease of $L$ - since $m \ll qT$, so that $L=r_2CT/(qT+m)\approx r_2 C/ q$ -, in order for the constraints in Eqs.~\eqref{eq:Con1} and \eqref{eq:Con2} to hold.~In contrast, a potential increase in $a$ will have exactly the opposite effects: reducing the duration of the explosive stage, increasing $T$, decreasing $N$ and accelerating their rates of change.~These predictions, which were reached on the basis of the CSP diagnostics, are in perfect agreement with the response of the system in a 20\% reduction and increase of parameter $a$, as displayed in Fig.~\ref{fig:1_abP}a.

Furthermore, despite the fact the process 7 is a process that can potentially affect the evolution of $T$ as shown in Eq.~\eqref{eq:VecF}, the CSP diagnostics indicate that process 7 is expected to have negligible effects to the constraints that bound the evolution of the system, as well as to its slow dynamics; see the discussion of Tables~\ref{tb:M1diags} and \ref{tb:M3diags}.~Indeed, the prediction provided by the CSP tools is in perfect agreement with the response of the system in a 40\% reduction of $c$, as displayed in Fig.~\ref{fig:1_abP}b, according to which no particular effect is reported.

Finally, the introduction of a potential perturbation on the rate constant of process 2 is expected to have no effect at the slow dynamics of the system, since process 2 is not driving the slow system; see the discussion of Table~\ref{tb:M1diags}.~However, it is expected to affect the constraint in Eq.~\eqref{eq:Con1}, which implies $N\approx eC/(pT)$.~Thus, a potential reduction in $e$ is expected to reduce $N$ (the mostly related to the 2nd constraint variable), in order for the constraint to hold.~Again, such a behavior is in perfect agreement with the response of the system in a 40\% reduction of $e$, as displayed in Fig.~\ref{fig:1_abP}b.~Indeed, no effect is noted to the duration of the explosive stage, or to the slow evolution of $T$, but only at the evolution of $N$, that follows a similar bounded - to the constraints - behavior with the reference case, at lower population levels.

%%%%%%%%%%%%%%%%%%%%%%%%%%%%%%%%%%%%%%%%
\section{The reduced model at the explosive stage}
\label{sec:redmod}

After identifying in Section~\ref{sec:CSPdyn} the constraints on which the system is bound to evolve and the processes governing its slow evolution along them, here a reduced model is constructed algorithmically for the explosive stage of the tumor-immune system dynamics.

Due to the fast/slow timescale separation during the explosive stage, the detailed model in Eq.~\eqref{eq:VecF} can be reduced to a reduced model of the form of Eq.~\eqref{eq:gov3}.~The 2-dim. SIM is defined by the algebraic equations of the constraints in Eqs.~\eqref{eq:Con1} and \eqref{eq:Con2}, which were shown to accurately approximate the fast dynamics and the related $N$ and $L$ cell populations in different initial conditions and parameter values; see Figs.~\ref{fig:TP_Sol+SIM_err9} and S1.~Despite having examined the governing slow dynamics, here the \textit{Importance Index} (II) in Eq.~\eqref{eq:gov8} is utilized in order to identify the processes contributing to the slow production/consumption of each cell population, so that the differential equations of the slow system in Eq.~\eqref{eq:gov3} can be defined.

The processes identified by II to contribute to the slow production/consumption of $T$ and $C$ populations and of $N$ and $L$ populations during the explosive stage are displayed in Figs.~S2 and S3 respectively for all the cases considered.~A summary of the most significant processes ($II_k^i>2\%$) is displayed in Table~\ref{tb:II}, in which the processes contributing after the initial transient phase and up to the very end of the explosive stage $t/t_{exp}<0.95$ are denoted in bold.

\begin{table}[!t]
\scriptsize
\centering
\begin{tabular}{c l}  \toprule[1.5pt]
Cell population   & Slow System  \\[2pt] \midrule[1pt]
$\dfrac{dT}{dt}=$ & $f(+\mathbf{R^1},-\mathbf{R^8},-(\mathbf{R^{12}}-\mathbf{R^{14}}))$   \\[6pt]
$\dfrac{dC}{dt}=$ & $f(+\mathbf{R^3},-\mathbf{R^6})$   \\[6pt]
$\dfrac{dN}{dt}=$ & $f(-\mathbf{R^1},+\mathbf{R^3},-\mathbf{R^6},+\mathbf{R^8},+(\mathbf{R^{12}}-\mathbf{R^{14}}),+(R^2-R^{13}))$   \\[6pt] 
$\dfrac{dL}{dt}=$ & $f(+\mathbf{R^1},+\mathbf{R^3},-R^5,-\mathbf{R^6},-\mathbf{R^8},+(\mathbf{R^{12}}-\mathbf{R^{14}}),-R^{15})$   \\[6pt]
\bottomrule[1.5pt]
\end{tabular}
\caption{Summary of the processes identified by II to contribute for the production (positive) or the consumption (negative) of each cell population for all the TP, TR, TP1 and TR1 cases considered.~The processes contributing after the initial transient phase and at the major part of the explosive stage are denoted in bold; $0.1<t/t_{exp}<0.95$.}
\label{tb:II}
\end{table}

Since $T$ and $C$ are not related to the formation of the constraints in Eqs.~\eqref{eq:Con1} and \eqref{eq:Con2}, their evolution during the explosive stage is to leading order governed by the slow system in Table~\ref{tb:II}, according to which:
\begin{itemize}
\item the tumor cell population $T$ is mainly produced by tumor growth (process 1) and consumed by the fractional cell kill rate by CD8$^+$ T cells (process 8).~Higher order corrections tending to consume $T$ are provided by the net rate $R^{12}-R^{14}$ of processes 12 and 14, which equilibrate for the formation of the 2nd constraint in Eq.~\eqref{eq:Con2}.
\item the circulating lymphocytes cell population $C$ evolves according to their growth and death rates (processes 3 and 6).
\end{itemize}
The evolution of $N$ and $L$ during the explosive stage is primarily dictated by the constraints in Eqs.~\eqref{eq:Con1} and \eqref{eq:Con2}.~However, higher order corrections are provided by the processes contributing to their slow evolution, as shown in Table~\ref{tb:II}, according to which:
\begin{itemize}
\item the NK cell population $N$ is mainly produced by processes 3, 8 and the net rate of processes 12 and 14 and consumed by processes 1 and 6.~Additional contributions are provided by the net rate of processes 2 and 13, but only during the initial transient phase and at the very end of the explosive stage.
\item the CD8$^+$ T cell population $L$ is mainly produced by processes 1, 3  and the net rate of processes 12 and 14, $R^{12}-R^{14}$ and consumed by processes 6 and 8.~Additional contributions are provided by processes 5 and 15, but only at the very end of the explosive stage.
\end{itemize}

The above identifications are in perfect agreement to the governing dynamics at the explosive stage, as identified by the TPI of the explosive CSP mode in Table~\ref{tb:M3diags}.~In particular, the driving processes 1, 8 and the net rate of processes 12 and 14 was identified in Table~\ref{tb:II} to play a major role on the slow system of all the cell populations - except from $C$ that is the slowest cell population, neither being related to the constraints, nor to the characteristic explosive CSP mode.

Finally, note that towards the end of the explosive stage, the role of process 1 is more significant at the tumor persistence TP and TP1 cases, in contrast to the tumor remission TR and TR1 cases, in which the role of process 8 is more decisive, according to Fig.~S1.~This feature indicates that the slow dynamics of tumor progression is governed by tumor growth in tumor persistence cases, in contrast to tumor remission ones where the tumor cell kill rate by CD8$^+$ T cells dominates the slow dynamics; this behavior was also reported by the TPI of the explosive CSP mode in Table~\ref{tb:M3diags}.

In summary, according to the CSP identifications of (i) the constraints in Eq.~\eqref{eq:Con1} and \eqref{eq:Con2} and (ii) the processes contributing to the slow production/consumption of the cell populations in Table~\ref{tb:II}, the reduced model that is valid with higher order accuracy during the explosive stage is:
\begin{align}
N & =  \dfrac{e C}{ p T}     &   \dfrac{dN}{dt} & =  [-a T (1 - b T) +  DT + \alpha - \beta C + (r_2 C T - q L T)]    \nonumber \\
L & =  \dfrac{r_2 C T}{q T+m} &  \dfrac{dL}{dt} & =  [a T (1 - b T) -  DT + \alpha - \beta C + (r_2 C T - q L T)]  \nonumber \\
&	&	\dfrac{dT}{dT} & = a T (1 - b T) - DT -  [(r_2 C T - q L T)] \nonumber \\
&	&	\dfrac{dC}{dT} & = \alpha - \beta C \label{eq:RedMod_HOA}
\end{align}
where the terms in brackets $[\cdot]$ provide higher order accuracy.~However, limited to leading order accuracy, the reduced model of Eq.~\eqref{eq:RedMod_HOA} can be further simplified to:
\begin{equation}
 N = \dfrac{e C}{ p T}, \qquad  L = \dfrac{r_2 C T}{q T+m}, \qquad \dfrac{dT}{dt}= a T (1 - b T) - DT, \qquad \dfrac{dC}{dt} = \alpha - \beta C 
\label{eq:RedMod_LOA}
\end{equation}
where $D=d \left( (L/T)^l \right) / \left( s+(L/T)^l \right)$.~The reduced model in Eq.~\eqref{eq:RedMod_LOA} consists of 2 algebraic equations for the cell populations related to the fast dynamics, $N$ and $L$, and 2 differential equations for the slow system formulating the governing slow evolution of $T$ and $C$.~The model essentially includes 10 parameters: $e/p$, $q/r_2$, $m/r_2$ for the algebraic equations and 7 for the differential ones and reproduces to leading order all the insights of the detailed model during the explosive stage.

%%%%%%%%%%%%%%%%%%%%%%%%%%%%%%%%%%%%%%%%%%%%%%%%%%%%%%%%%%%%%%%%%%%%%%%%%%%%%%%%
\section{Conclusions}
\label{sec:Con}

In this work, the interactions between the tumor and the immune system were investigated on the basis of a fundamental mechanistic model \cite{de2006mixed}, which comprehends the basic dynamical features of many cancer immunology models: complex interactions; multi-scale dynamics; overparameterized character.~In order to tackle these features, which traditionally hinder the accessibility and wide utilization of models in mathematical biology, a systematic and algorithmic framework in the context of GSPT was introduced with the aim to (i) provide mechanistic understanding, (ii) reveal the underlying mechanisms controlling tumor progression and (iii) simplify complexity by delivering reduced models.

The model's equilibria and their stability were investigated by validating that the model successfully reproduces the long-term tumor behavior; i.e., the tumor escaping immune surveillance (HTE) and the tumor being suppressed by the immune system (TFE).~Aiming to reveal the processes driving tumor progression in Section~\ref{sec:exp}, the existence of an initial \emph{explosive stage} was discovered, during which explosive timescales characterize the evolution of the system.~The explosive stage is crucial for determining the basins of attraction of the TFE and HTE, since the unstable HTE threshold branch is not enough to \textit{a priory} determine tumor persistence or remission, as highlighted by the cases considered in Fig.~\ref{fig:RPSolComp}.~More importantly, such a result indicates that the long-term evolution of the tumor is determined during the explosive stage of tumor progression.~Note that the crucial role of the explosive stage to the system's dynamics has been identified in a number of studies mostly related to reactive flows \citep{manias2016mechanism,valorani2006automatic,tingas2018chemical}, but also to COVID-19 pandemics \citep{patsatzis2021relation} and oscillators' and attractors' dynamics \citep{maris2015hidden,goussis2013role}; this is the first time reported to cancer immunology dynamics.

Exploiting the multi-scale character of the model, the CSP algorithmic tools were employed in Section~\ref{subsec:fast} in order to identify the constraints along which the system is bound to evolve during the explosive stage.~It was shown that the underlying mechanism dictating tumor progression evolves along two constraints: the 1st related to NK cell population and formed between their growth and inactivation processes and the 2nd related to CD8$^+$ T cell population and formed mainly between their recruitment and inactivation processes; see Table~\ref{tb:M1diags}.~The validity of the constraints was further examined in a number of tumor persistence/remission cases and in different parameter sets corresponding to patient 9 and 10 \cite{de2006mixed}, indicating the preservation of the confining tumor progression dynamics along different initial conditions and parameter values. 

More importantly, it was identified by CSP and further cross-validated by indicative perturbations in Section~\ref{subsec:slow} that the slow dynamics of tumor progression during the explosive stage (i) is majorly governed by the tumor growth rate and to a much lesser degree by the fractional cell kill rate of tumor cells by NK cells and (ii) primarily affects the slow evolution of tumor cell population.~In addition, by examining a number of tumor persistence/remission cases at different parameter sets, it was further identified that the recruitment and inactivation of CD8$^+$ T cells provide additional contributions to the governing slow dynamics, which are more pronounced at the tumor remission cases.

On top of these insights, it was further demonstrated the ability of the CSP algorithmic framework to provide accessible reduced models in Section~\ref{sec:redmod}.~In particular, the CSP-constructed reduced model in Eq.~\eqref{eq:RedMod_LOA} requires 10 fewer parameters to accurately (by leading order) provide the insights of the detailed model in Eq.~\eqref{eq:VecF} during the critical explosive stage of tumor progression.~Constituting by 2 algebraic and 2 differential equations, rather than 4 differential ones and being free of fast timescales \citep{goussis2012quasi}, the reduced model in Eq.~\eqref{eq:RedMod_LOA} provides significant simplification that could not otherwise reached; employing the traditional QSSA for $N$ and $L$ would result in constructing a much complexer reduced model, since the processes contributing to a negligible degree (such as process 7, 9, 10, 11) would be incorporated in the model.~Last but not least, the reduced model in Eq.~\eqref{eq:RedMod_LOA} is a simplified version of the ones constructed by the algorithmic procedure provided by CSP, which guarantees their validity and accuracy under specific algorithmic criteria \citep{goussis2012quasi,patsatzis2019new}.

Finally, by successfully extending the application of CSP to cancer immunology models, the framework employed in this work provides future perspectives mainly towards three directions.~First, due to its algorithmic nature, it is not hindered by the size and the complexity of the model under consideration, thus rendering its application straightforward in more detailed models.~Secondly, it was demonstrated in Section~\ref{subsubsec:demo} that the CSP identifications can successfully predict the impact of parameter variations to the model's outcome, accounting for the mutli-scale dynamics of the model, in contrast to sensitivity analysis that neglects it \citep{li16,hu2018dynamics,wilson12,rihan14}.~The insight provided on the processes and cell types that majorly affect the model's outcome is particularly significant when accounting for anti-cancer therapies, since they can indicate which interactions and cell types to target for more effective treatment regimes.~Last but not least, the algorithmic model reduction employed in this work can successfully provide accessible reduced models of fewer parameters and free of fast timescales, which are valid over a wider parameter range than the ones provided by the traditional paper-pencil techniques \citep{goussis2012quasi,patsatzis2019new}.~Since the latter techniques are employed by hand, in cancer immunology as well \citep{osojnik20,kuznetsov94,louzoun2014mathematical}, the algorithmic nature of this framework can serve as the final step of modeling pipeline for providing valid simplified models, so that the development of detailed models is no longer hindered by overparameterization.

\appendix
% equation numbering
%\numberwithin{equation}{section}
%\setcounter{equation}{0}
% figure numbering
\numberwithin{figure}{section}
\setcounter{figure}{0}
% table numbering
\numberwithin{table}{section}
\setcounter{table}{0}

%%%%%%%%%%%%%%%%%%%%%%%%%%%%%%%%%%%%%%%%%%%%%%
\section{The equilibria of the system and their stability}
\label{app:Eq+St}

Let the state variables of cell populations at equilibrium be $\mathbf{y}^*=(T^*, N^*, L^*, C^*)$.~Setting the differential equations of $N$ and $C$ in Eq.~\eqref{eq:VecF} to zero, implies:
\begin{equation}
C^* =\dfrac{\alpha}{\beta} \qquad N^*=\dfrac{\alpha e \left( h+T^{*2}\right)}{\beta \left( fh + hpT^* + (f-g) T^{*2} + pT^{*3}\right)}\equiv N^*(T^*)
\label{eq:SP1}
\end{equation}
so that $C^*$ is constant and $N^*$ is a function of $T^*$.~Similarly, setting the differential equations of $T$ and $L$ in Eq.~\eqref{eq:VecF} to zero allows for the determination of $T^*$ and $L^*$, which requires the solution of the 2-dim. system of algebraic equations: 
\begin{align}
a T^* (1-bT^*) - cN^*T^* - D^*T^*  & = 0 \label{eq:SP2} \\
\alpha^*(T^*) L^{*2} +  \beta^*(T^*) L^* + \gamma^*(T^*) &=0 \label{eq:SP3}
\end{align}
where $D^*=d( L^*/T^*)^l/\left(s+\left( L^*/T^*\right)^l\right)$ and
\begin{equation*}
\alpha^*(T^*) \equiv  -uN^* \qquad \beta^*(T^*) \equiv   -m + j\dfrac{(D^*T^*)^2}{k+(D^*T^*)^2} -qT^* \qquad  \gamma^*(T^*) \equiv  (r_1 N^* +r_2 C^*)T^*
\end{equation*}
Note that although $D^*$ is a function of both $T^*$ and $L^*$, Eq.~\eqref{eq:SP2} implies that $D^* T^*$ is function of $T^*$ only, since $N^*=N^*(T^*)$ according to the second relation in Eq.~\eqref{eq:SP1}.~As a result, the coefficients of the quadratic Eq.~\eqref{eq:SP3}, $\alpha^*$, $\beta^*$ and $\gamma^*$ are also functions of $T^*$.~The non-linear system of Eqs.~\eqref{eq:SP2} and \eqref{eq:SP3} has multiple real and complex solutions that can be obtained numerically in general.

In order to examine the existence of tumor-free equilibria, it is firstly noted that $T=0$ is precluded from the domain of the system in Eq.~\eqref{eq:VecF}, due to the definition of $D$.~However, the limiting case where $T$ approaches to zero asymptotically ($T\rightarrow 0$) is permitted and expresses the absence of tumour, since $T$ represents the number of tumour cells.~In the limit $T^*\rightarrow0$, Eq.~\eqref{eq:SP2} is satisfied and Eq.~\eqref{eq:SP3} is simplified to:
\begin{equation}
-u N^* L^{*2}  -m L^*  =0 \Rightarrow L^* = 0 \quad \text{or} \quad L^*=-m/ (u N^*)
\end{equation} 
Substituting from Eq.~\eqref{eq:SP1}, two tumor-free equilibria arise: the biologically feasible $\mathbf{y}^{*1}=(0,\alpha e/\beta f,0,\alpha/\beta)$ and the biologically infeasible $\mathbf{y}^{*2}=(0,\alpha e/\beta f,-(m f \beta)/(u e \alpha),\alpha/\beta)$, since $L^*<0$.~At the tumor-free equilibrium $\mathbf{y}^{*1}$, the eigenvalues of the Jacobian matrix can be analytically calculated:
\begin{equation*}
\lambda_1 = -\dfrac{\alpha c e}{\beta f} + a - d \qquad \lambda_2 = -f \qquad \lambda_3 = -m \qquad \lambda_4 = - \beta
\end{equation*}
Since all the parameters are positive and $\lambda_2$, $\lambda_3$ and $\lambda_4$ are always negative, the tumor-free equilibrium $\mathbf{y}^{*1}$ is stable if and only if $\lambda_1<0 \Leftrightarrow (a-d) \beta f< \alpha c e$.~On the other hand, the unfeasible tumor-free equilibrium $\mathbf{y}^{*2}$ is always unstable, since the eigenvalue $\lambda_3=m$ is always positive.

Turning the focus in high-tumor equilibria, the conditions $T^*,N^*,L^*> 0$ were imposed.~Within these limits, the discriminant of the quadratic Eq.~\eqref{eq:SP3} is always positive, so that two distinct real solutions $L_{\pm}^*$ arise; in particular, $L_+^*<0$ and $L_-^*>0$, since $\gamma^*(T^*)/\alpha^*(T^*)$ is always negative.~Since $L_-^*$ is the only biological feasible choice, its substitution in Eqs.~\eqref{eq:SP2} and \eqref{eq:SP1} for the numeric calculation of $T^*>0$ and $N^*>0$, respectively, provides all the biologically feasible high-tumor equilibria of the system.~The stability of these equilibria is provided by examining the eigenvalues of the related Jacobian matrices, which are also calculated numerically.

%%%%%%%%%%%%%%%%%%%%%%%%%%%%%%%%%%%%%%%%%%%%%%
\section{The parameter set of the mathematical model}
\label{app:ParSet}

The analysis of the mathematical model in Eq.~\eqref{eq:VecF}, which incorporates the processes enlisted in Table~\ref{tb:RR}, was carried out by adopting the parameter set of patient 9 in Ref.~\cite{de2006mixed}, as enlisted in Table~\ref{tb:par9}.

\begin{table}[!h]
\centering
\begin{tabular}{c l c}
\multicolumn{2}{c}{Parameter value} & Unit \\[2pt]
\hline
$a$ 	& 4.31 $\times$ 10$^{-1}$		& $1/day$ \\
$b$ 	& 1.02 $\times$ 10$^{-9}$ 	& $1/cell$ \\
$e$ 	& 2.08 $\times$ 10$^{-7}$ 	& $1/day$ \\
$\alpha$	& 7.50 $\times$ 10$^{+8}$	& $cell/day$ \\
$f$  	& 4.12 $\times$ 10$^{-2}$ 	& $1/day$ \\
$m$ 	& 2.04 $\times$ 10$^{-1}$		& $1/day$ \\
$\beta$	& 1.20 $\times$ 10$^{-2}$		& $1/day$ \\
$c$	& 6.41 $\times$ 10$^{-11}$		& $1/(day\ cell)$\\
$d$	& 2.34 		& $1/day$\\
$l$	& 2.09 		& -
\end{tabular}
\ 
\begin{tabular}{c l c}
\multicolumn{2}{c}{Parameter value} & Unit \\[2pt]
\hline
$s$	& 8.39 $\times$ 10$^{-2}$ 		& - \\
$g$ 	& 1.25 $\times$ 10$^{-2}$		& $1/day$ \\
$h$ 	& 2.02 $\times$ 10$^{+7}$		& $cell^2$ \\
$j$ 	& 2.49 $\times$ 10$^{-2}$		& $1/day$ \\
$k$ 	& 3.66 $\times$ 10$^{+7}$		& $cell^2$ \\
$r_1$ 	& 1.10 $\times$ 10$^{-7}$		& $1/(day\ cell)$ \\
$r_2$ 	& 6.50 $\times$ 10$^{-11}$		& $1/(day\ cell)$ \\
$p$ 	& 3.42 $\times$ 10$^{-6}$		& $1/(day\ cell)$ \\
$q$ 	& 1.42 $\times$ 10$^{-6}$		& $1/(day\ cell)$ \\
$u$ 	& 3.00 $\times$ 10$^{-10}$		& $1/(day\ cell^2)$ 	
\end{tabular}
\caption{The parameter set of patient 9 in Ref.~\cite{de2006mixed}, considered for the analysis of the model in Eq.~\eqref{eq:VecF}}
\label{tb:par9}
\end{table}

%%%%%%%%%%%%%%%%%%%%%%%%%%%%%%%%%%%%%%%%%%%%%%
\section{The CSP algorithmic methodology and its diagnostic tools}
\label{app:CSP}

The \emph{Computational Singular Perturbation} (CSP) methodology provides a systematic and algorithmic framework to deliver asymptotic analysis.~It is employed for the analysis of multi-scale dynamical systems, in which the various processes incorporated in the model, act in a wide range of timescales.~CSP exploits the fast/slow separation of timescales by algorithmically decomposing the tangent space, in which the system evolves, in a fast and a slow subspace; thus identifying the components of the system (processes and variables) that essentially generate the fast and the slow dynamics.~Such identifications are of great significance to determine the trending dynamics of the system, since the fast components tend to drive the system towards an equilibrium.~This feature is manifested by the emergence of specific constraints, along which the system is confined to evolve, governed by the action of the slow components.~Thus, the decomposition in a fast and a slow subspace provided by CSP, enables the determination of the processes that contribute to the emergence of the constraints, as well as the processes driving the slow dynamics.~Such an algorithmic determination is particularly insightful in systems incorporating many processes, especially when the optimal goal is to control the system dynamics.

Considering a $N$-dim. dynamical system of Ordinary Differential Equations (ODE), such as the one in Eq.~\eqref{eq:VecF}, the first step is to set it in its vector form:
\begin{equation}
\dfrac{d\mathbf{y}}{dt}  = \mathbf{g}(\mathbf{y}) = \sum_{k=1}^{K}  \mathbf{S}_k R^k(\mathbf{y})
\label{eq:gov1}
\end{equation}
where $\mathbf{y}$ is the $N$-dim.~\emph{state vector} containing the variables and $\mathbf{g(y)}$ is the $N$-dim.~vector field consisting of $K$ processes, the stoichiometric vectors of which are $\mathbf{S}_k$ and the related process rates  are $R^k(\mathbf{y})$ for $k=1,\dots,K$.~CSP decomposes the vector field $\mathbf{g(y)}$ in $N$ modes \citep{lam1989,lam1991conventional}:
\begin{equation}
\dfrac{d\mathbf{y}}{dt}  =  \mathbf{g}(\mathbf{y}) =  \sum_{n=1}^{N}\mathbf{a}_n (\mathbf{y})f^n(\mathbf{y}), \qquad 
f^n (\mathbf{y})  = \mathbf{b}^n(\mathbf{y}) \cdot \mathbf{g(y)} = \sum_{k=1}^{K} \left( \mathbf{b}^n (\mathbf{y}) \cdot \mathbf{S}_k\right)R^k(\mathbf{y}) \label{eq:gov2} 
\end{equation} 
by introducing the $N$-dim.~CSP column basis vectors $\mathbf{a}_n (\mathbf{y})$ of the $n$-th mode and their $N$-dim.~dual row vectors $\mathbf{b}^n(\mathbf{y})$; the latter satisfying the orthogonality conditions $\mathbf{b}^i(\mathbf{y}) \cdot \mathbf{a}_j (\mathbf{y}) = \delta_j^i$ \citep{lam1989,lam1994}.~In this way, each CSP mode $\mathbf{a}_n(\mathbf{y})f^n(\mathbf{y})$ in Eq.~\eqref{eq:gov2} is characterized by a timescale $\tau_n(\mathbf{y})$ and an amplitude $f^n(\mathbf{y})$: the timescale $\tau_n(\mathbf{y})$ providing a measure of the time frame of action of the $n$-th CSP mode and the amplitude $f^n(\mathbf{y})$ providing a measure of the impact of the mode in driving the trajectory along the direction of $\mathbf{a}_n(\mathbf{y})$.

Due to the fast/slow timescale separation, the system exhibits, say $M$, timescales that are (i) dissipative in nature; i.e., their action tends to drive system towards to the stable equilibrium and (ii) much faster than the rest.~In this case, the CSP form in Eq.~\eqref{eq:gov2} can be reduced to the system of Differential Algebraic Equations (DAE):
\begin{equation}
f^r(\mathbf{y}) \approx 0\quad(r=1,\ldots,M), \qquad  
\dfrac{d\mathbf{y}}{dt} \approx \sum_{s=M+1}^{N}\mathbf{a}_s(\mathbf{y})f^s(\mathbf{y}) 
\label{eq:gov3}
\end{equation}
that consists the reduced model, which is valid when the $M$ timescales are exhausted; i.e., when the $M$ constraints have emerged, so that they confine the solution to evolve along them.~The $M$-dim.~system of algebraic equations in Eq.~\eqref{eq:gov3} defines the \emph{Slow Invariant Manifold} (SIM), that is a low dimensional surface emerging in phase-space, where the system is confined to evolve \citep{fenichel1979geometric,kaper1999systems,verhulst2006nonlinear,hek2010geometric}.~The $N$-dim.~system of differential equations in Eq.~\eqref{eq:gov3} defines the \emph{slow system} that governs the slow dynamics on the SIM.~Note that when the reduced model in Eq.~\eqref{eq:gov3} is valid, the slow system is free of fast timescales; thus reproducing the slow dynamics incorporated in the full ODE model in Eq.~\eqref{eq:gov1} \citep{goussis2012quasi,patsatzis2019new,maris2015hidden}.

The decomposition of the tangent space to a fast and a slow subspace provided by CSP requires the calculation of the CSP basis vectors.~CSP provides an algorithmic procedure for the calculation of the CSP basis vectors via two iterative refinements \citep{lam1989,lam93,goussis99}.~However, a leading order accurate estimation of the CSP vectors $\mathbf{a}_i(\mathbf{y})$ and $\mathbf{b}^i(\mathbf{y})$ ($i=1,\ldots,N$) is provided by the right $\boldsymbol{\alpha}_i(\mathbf{y})$ and left $\boldsymbol{\beta}^i(\mathbf{y})$, respectively, eigenvectors of the $N \times N$-dim.~Jacobian $\mathbf{J}(\mathbf{y})$ of $\mathbf{g}(\mathbf{y})$ \citep{lam1991conventional,lam1989,tingas2016ignition}.~In the following, the CSP diagnostic tools are presented considering $\mathbf{a}_i=\boldsymbol{\alpha}_i$ and $\mathbf{b}^i=\boldsymbol{\beta}^i$ and the dependency from $\mathbf{y}$ is dropped for simplicity.

\subsection{The CSP diagnostic tools}
\label{appSub:tools}

CSP provides a number of diagnostic tools that are utilized to acquire the relevant physical understanding of the system under consideration.~In this work, the CSP tools were employed to identify (i) the physical processes contributing to the formation of the emerging constraints, (ii) the processes generating the fast/slow timescales, (iii) the cell populations related to the development of fast/slow timescales and (iv) the processes contributing to the production/consumption of each variable along the constraints; i.e., the governing processes of the slow system.

The $M$ exhausted modes impose the emergence of the $M$ constraints, which originate as a result of significant cancellations among the various processes.~Each amplitude in the first expression of Eq.~\eqref{eq:gov3} can be written in the form:
\begin{equation}
f^r = (\boldsymbol{\beta}^r \cdot \mathbf{S}_1)R^1 + \ldots +  (\boldsymbol{\beta}^r  \cdot \mathbf{S}_K)R^K \approx 0
\label{eq:gov4}
\end{equation}
for $r=1,\ldots,M$, where the analytic expression of the amplitudes in Eq.~\eqref{eq:gov2} is utilized.~Some of the additive terms in Eq.~\eqref{eq:gov4} introduce non-negligible contributions, some of which cancel each other, thus forming the emerging constraints.~In order to identify these terms, the relative contribution of the $k$-th process ($k=1,\ldots,K$) to each fast amplitude $f^r\approx 0$ ($r=1,\ldots,M$) is measured by the \emph{Amplitude Participation Index} (API):
\begin{equation}
P^r_k = \dfrac{ (\boldsymbol{\beta}^{r} \mathbf{S}_k)R^k}{ \sum\nolimits_{i=1}^{K}|    (\boldsymbol{\beta}^{r}\mathbf{S}_i)R^i         | }
\label{eq:gov5}
\end{equation}
where by definition $\sum_{k=1}^{K}|P_k^{r}|=1$,~\citep{lam1994,goussis2006b,valorani2003}.~$P^r_k$ can be either positive or negative and the sum of the positive and negative terms equals by definition to 0.5 and -0.5, respectively.

The time frame in which each of the $M$ constraints is formed is characterized by the related to this constraint fast timescale.~In addition, the evolution of the system along the constraints is characterized by the fastest of the slow $N-M$ timescales, frequently called the \emph{characteristic} timescale, $\tau_{M+1}$.~Both the fast and the slow timescales are estimated by the inverse norm of the eigenvalues of the Jacobian $\mathbf{J}$ of $\mathbf{g}$; i.e., $\tau_n = |\lambda_n|^{-1}$ for $n=1,\ldots,N$.~In order to identify the processes generating the fast/slow timescales, the relative contribution of the $k$-th process ($k=1,\ldots,K$) to the timescale $\tau_n$ is measured by the \emph{Timescale Participation Index} (TPI): 
\begin{equation}
J^n_k = \dfrac{c_k^n}{ \sum\nolimits_{i=1}^{K}|    c^n_i         |} \qquad \text{where} \qquad  c_k^n=\boldsymbol{\beta}^n \nabla  \left(  \mathbf{S}_kR^k \right)  \boldsymbol{\alpha}_n 
\label{eq:gov6}
\end{equation}
where by definition $\sum_{k=1}^{K}|J_k^n|=1$~\citep{goussis2005a,goussis2006b,diamantis2015h2} and the term $c_k^n$ expresses the contribution of the $k$-th process to the $n$-th eigenvalue ${\lambda}_n =  c_1^n + \ldots + c_{K}^n$, since $\mathbf{J}=\sum_{k=1}^K \nabla  \left(  \mathbf{S}_kR^k \right)$.~$J^n_k$ can be either positive or negative, implying that the $k$-th process  contributes to an explosive or dissipative character of the $n$-th timescale $\tau_n$.~When $\tau_n$ is explosive in nature, the processes with positive $J^n_k$ promote the explosive nature, while the ones with negative $J^n_k$ oppose to it (and vice versa for the dissipative timescales).~By definition, the dissipative (explosive) timescales relate to the components of the system that tend to drive it towards (away from) a stable equilibrium \citep{lam1989,lam1994}, since the character (dissipative/explosive) of a timescale is determined by the real part of the respective eigenvalue (negative/positive).

Each variable associates differently to each exhausted CSP mode and thus the related timescale; e.g., the variables considered as fast relate mostly to a fast CSP mode ($m=1, \dots, M$) and to a much lesser degree with the remaining slow CSP modes.~The relation of each variable to the  $m$-th CSP mode is identified by the \emph{CSP Pointer} (Po):
\begin{equation}
\mathbf{D}^m = diag \left[ \boldsymbol{\alpha}_m\boldsymbol{\beta}^m \right] = \left[ \alpha^1_m\beta^m_1,\alpha^2_m\beta^m_2,\ldots,\alpha^{N}_m\beta^m_{N} \right] 
\label{eq:gov7}
\end{equation}
where $\sum_{i=1}^N\alpha^i_m\beta^m_i=1$, due to the orthogonality condition $\boldsymbol{\beta}^i \cdot \boldsymbol{\alpha}_j=\delta^i_j$~\citep{goussis1992study,lam1994,valorani2003,goussis2012quasi}.~Large values of the $i$-th element of Po, $\alpha^i_m\beta^m_i$, imply strong correlation of the $i$-th variable to the $m$-th CSP mode, while a value close to unity implies that the $i$-th variable is potentially in \emph{Quasi Steady-State} (QSS) \citep{goussis2012quasi}.

Finally, the processes contributing the most to the slow evolution of each variable along the $M$ constraints, according to the $N$-dim. slow system of ODEs in Eq.~\eqref{eq:gov3}, are identified by the slow \textit{Importance Index} (II):
\begin{equation}
II_k^n=\dfrac{\sum\nolimits_{s=M+1}^{N}\alpha_s^n (\boldsymbol{\beta}^s\cdot\mathbf{S}_k)R^k}{\sum\nolimits_{j=1}^K|\sum\nolimits_{s=M+1}^{N}\alpha_s^n (\boldsymbol{\beta}^s\cdot\mathbf{S}_j)R^j|}
\label{eq:gov8}
\end{equation}
where by definition $\sum\nolimits_{k=1}^K|II^n_k|=1$ for each of the $n=1,\ldots,N$ variables \citep{lam1994,goussis2006b}.~$II_k^n$ measures the relative contribution of the $k$-th process to the evolution of the $n$-th variable according to the slow system in Eq.~\eqref{eq:gov3}.~$II_k^n$ can be either positive or negative, implying the importance of the $k$-th reaction to the production or consumption of the $n$-th variable \citep{lam1994,valorani2003,goussis1992study,maris2015hidden}.

%%%%%%%%%%%%%%%%%%%%%%%%%%%%%%%%%%%%%%%%%%%%%%%%%%%%%%%%%%%%%%%%%%%%%%%%%%%%%%%%
%\section*{References}

\bibliography{mybibfile}

\end{document}